\documentclass[amsfonts,12pt]{article}

\usepackage{amsfonts,amsmath}
\usepackage{hyperref}

\newtheorem{theorem}{Theorem}

\newtheorem{lemma}[theorem]{Lemma}
\newtheorem{proposition}[theorem]{Proposition}
\newtheorem{definition}[theorem]{Definition}
\newtheorem{example}[theorem]{Example}

\newtheorem{remark}{Remark}
\DeclareMathOperator{\Res}{Res}

\def\ppp{{\mathbb{P}}}

\def\fff{{\mathbb{F}}}
\def\eee{{\mathbb{E}}}
\def\qqq{\mathbb{Q}}

\def\ccc{\mathbb{C}}
\def\zzz{\mathbb{Z}}

\def\pf{{\bf proof}:\ }
\def\qed{$\Box$}

\begin{document}

\author{David Joyner and Amy Ksir\thanks{Mathematics Dept, USNA, Annapolis, MD 21402,
wdj@usna.edu and ksir@usna.edu}
\thanks{This is the expanded text of a talk given at the
Algorithmic Algebraic Geometry session of the
AMS Meeting in Altanta, GA, January 4-8, 2005.
}}
\title{Modular representations on some Riemann-Roch spaces \\
of modular curves $X(N)$
}
\date{2-28-2005}

\maketitle

\begin{abstract}
We compute the $PSL(2,N)$-module structure of the Riemann-Roch space
$L(D)$, where $D$ is an invariant non-special divisor on the modular
curve $X(N)$, with $N\geq 7$ prime.  This depends on a computation
of the ramification module, which we give explicitly.  These results
hold for characteristic $p$ if $X(N)$ has good reduction mod $p$ and
$p$ does not divide the order of $PSL(2,N)$.  We give as examples
the cases $N=7, 11$, which were also computed using \cite{GAP}.
Applications to AG codes associated to this curve are considered,
and specific examples are computed using \cite{GAP} and
\cite{MAGMA}.

\end{abstract}

\vskip .5in

\tableofcontents

\vskip .3in

\section{Introduction}

The modular curve $X(N)$ has a natural action by the finite group
$G=PSL(2,N)$, with quotient $X(1)$. If $D$ is a $PSL(2,N)$-invariant
divisor on $X(N)$, then there is a natural representation of
$G$ on the Riemann-Roch space $L(D)$. In this paper, we
give some results about the $PSL(2,N)$-module structure of the
Riemann-Roch space $L(D)$, in the  case where $N$ is prime and $N
\geq 7$. If $D$ is non-special then a formula in
Borne \cite{B} gives

\begin{equation}
\label{eqn:Borne}
[L(D)]=(1-g_{X(1)})[k[G]]+[deg_{eq}(D)]-[\tilde{\Gamma}_G].
\end{equation}
Here $g_{X(1)}$ is the genus of $X(1)$ (which is zero), square
brackets denote the equivalence class of a representation of $G$,
$deg_{eq}(D)$ is the equivariant degree of $D$, and
$\tilde{\Gamma}_G$ is the ramification module  (these will be
defined in sections \ref{sec:rammod} and \ref{sec:equivdeg}).

Our main result (see Theorem \ref{thrm:main} in section
\ref{sec:rammod}) gives an explicit computation of the
$G$-module structure of the ramification module $\tilde{\Gamma}_G$.
We show it
can be explicitly decomposed into
irreducible $G$-modules where the multiplicity
of an irreducible $G$-module $\pi$ in
$\tilde{\Gamma}_G$ only
depends on $\pi$ and the residue class of $N$ modulo $24$.
We then go on in section \ref{sec:equivdeg} to compute the
equivariant degrees of $G$-invariant divisors on $X(N)$, and use
these and Borne's formula to compute the $G$-module structure of the
corresponding Riemann-Roch spaces.

If $K=\qqq(G)$ denotes the abelian extension of $\qqq$ generated by
the character values of $G$ then ${\cal G}=Gal(K/\qqq)$ acts on the
set of irreducible representations of $G$. We call this the {\it
Galois action}. In the case where $\tilde{\Gamma}_G$ is invariant
under the Galois action, the authors have given a somewhat simpler
formula for $\tilde{\Gamma}_G$ in \cite{JK}.  In section
\ref{sec:rammod} (see Theorem \ref{thrm:2}), we also prove that the
ramification module is $\mathcal{G}$-invariant if and only if
$N\equiv 1 \pmod 4$.

As a corollary, it is an easy exercise now to compute explicitly the
decomposition

\[
H^1(X(N),k) = H^0(X(N),\Omega^1)\oplus \overline{H^0(X(N),\Omega^1)}
=L(K)\oplus \overline{L(K)},
\]
into irreducible $G$-modules, where $K$ is a canonical divisor.
This was discussed in \cite{KP} (over $k=\ccc$)
and \cite{Sc} (over the finite field $k=GF(N)$).
Indeed, Schoen observes that
the multiplicities of the irreducible representations
occurring in $H^1(X(N),k)$ can be interpreted in terms of
dimension of cusp forms and number of cusps on $X(N)$.

In section \ref{sec:examples}
we look at the examples $N=7,11$, using \cite{GAP}
to do many of the computations.  In the last section, applications
to AG codes associated to this curve are considered (\cite{MAGMA}
was used to do some of these computations).

{\it Notation}: Throughout this paper, $N>5$ is a prime,
$\fff = GF(N)$ is the finite field with $N$ elements,
and $G=PSL(2,N)$.

{\it Acknowledgements}:
We thank D. Prasad and R. Guralnick for enlightening
correspondence, and in particular for the
references \cite{KP} and  \cite{Ja}.

\section{Modular curves}
\label{sec:klein}

Let $H$ denote the complex upper half-plane, let $H^* = H\cup
\qqq\cup \{\infty\}$, and recall that $SL(2,\qqq)$ acts on $H^*$ by
fractional linear transformations. Let $X(N)$ denote the modular
curve defined over $\qqq$ whose complex points are given by $\Gamma
(N)\backslash H^*$, where

\[
\Gamma(N)=\{ \left(
\begin{array}{cc}
a & b\\
c & d
\end{array}
\right)\in SL_2(\zzz)\ |\ a-1\equiv d-1\equiv b\equiv c\equiv 0
\pmod{N}\}.
\]

Throughout this paper, we will assume that $N$ is prime and $N>6$.
In this case, the genus of $X(N)$ is given by the formula

\[
g=1+\frac{(N-6)(N^2-1)}{24}.
\]
For example, $X(7)$ is genus $3$ and $X(11)$ is genus $26$.

Let $N$ be a prime and, for $j\in \zzz/N\zzz$, let $y_j$ be
variables satisfying

\begin{equation}
\label{eqn:klein} y_j+y_{-j}=0, \ \ \ y_{a+b}y_{a-b}y_{c+d}y_{c-d}+
y_{a+c}y_{a-c}y_{d+b}y_{d-b}+ y_{a+d}y_{a-d}y_{b+c}y_{b-c}=0,
\end{equation}
for all $a,b,c,d \in  \zzz/N\zzz$.  These are Klein's equations
for $X(N)$ (see Adler \cite{A} or Ritzenthaler \cite{R1}).

\begin{example}
When $N=7$, this reduces to
$y_1^3y_2-y_2^3y_3-y_3^3y_1=0$,
the famous Klein quartic.

When $N=11$, the $20$ equations which arise
reduce to the $10$ equations

\[
\begin{array}{r}
-y_1^2y_2y_3+y_2y_4y_5^2+y_3^2y_4y_5=0,\\
-y_1^3y_4+y_2y_4^3-y_3^3y_5=0,\\
-y_1y_3^3-y_1^3y_5+y_2^3y_4=0,\\
-y_1^2y_3y_4+y_1y_3y_5^2+y_2^2y_4y_5=0,\\
-y_1^2y_2y_5+y_1y_3y_4^2-y_2^2y_3y_5=0,\\
y_1^3y_2-y_3y_5^3-y_4^3y_5=0,\\
y_1y_5^3-y_2^3y_3+y_3^3y_4=0,\\
-y_1y_2^2y_4+y_1y_4^2y_5+y_2y_3^2y_5=0,\\
y_1y_2^3+y_2y_5^3-y_3y_4^3=0,\\
y_1y_2y_3^2+y_1y_4y_5^2-y_2y_3y_4^2=0.
\end{array}
\]
\end{example}

The curve $X(N)$ over a field $k$ parametrizes pairs of an elliptic
curve over $k$ and a subgroup of order $N$ of the group structure on
the curve.  This can be extended to fields of positive
characteristic, if $X(N)$ has good reduction.  Since Klein's
equations have integer coefficients they can also be extended to an
arbitrary field $k$. However, Velu \cite{V} (see also Ritzenthaler
\cite{R3}) has shown that $X(N)$ has good reduction over fields of
characteristic $p$ where $p$ does not divide $N$ (in our case,
$p\not= N$, since $N$ is itself assumed to be a prime).

Let

\[
G=PSL_2(\zzz/N\zzz)\cong \overline{\Gamma(1)}/\overline{\Gamma(N)},
\]
where the overline denotes the image in $PSL_2(\zzz)$. This group
acts on $X(N)$. (In characteristic $0$, see \cite{S}, in
characteristic $\ell>0$, see \cite{R1}.) When $N>2$ is prime,
$|G|=N(N^2-1)/2$.  The quotient of $X(N)$ by the action of
$PSL(2,N)$ is $X(1) \cong \ppp^1$.

\begin{definition}
\label{def:good} When $X$ has good reduction to a finite field $k$
and, in addition, the characteristic $\ell$ of $k$ does not divide
$|G|$, we say that $\ell$ is {\bf good}.
\end{definition}

If $k$ is a field of good characteristic, the automorphism group of
$X(N)$ is known to be $PSL(2,N)$ \cite{BCG}.

The action of $G=SL_2(\zzz/N\zzz)$ on the set of points of the
projective curve defined by Klein's equations is described in
\cite{R1} (see also \cite{A2}, \cite{R2}). The element $g=\left(
\begin{array}{cc}
a & b\\
c & d
\end{array}
\right)\in G$ sends $(y_j)_{j\in \zzz/N\zzz}\in X(N)$ to
$(\rho(g)y_j)_{j\in \zzz/N\zzz}\in X(N)$, where

\[
\rho(g)(y_j)=\sum_{t\in \zzz/N\zzz}
\zeta^{b(aj^2+2jtc)+t^2cd}y_{aj+tc},
\]
where $\zeta$ denotes a primitive $N$-th root of unity in $k$.

\begin{remark}
When the formulas for the special cases
$\rho\left(
\begin{array}{cc}
0 & 1\\
-1 & 0
\end{array}
\right)$, $\rho\left(
\begin{array}{cc}
1 & 1\\
0 & 1
\end{array}
\right)$, and
$\rho\left(
\begin{array}{cc}
a & 0\\
0 & a^{-1}
\end{array}
\right)$, are written down separately, the similarity with the
Weil representation for $SL_2(\zzz/N\zzz)$ is
striking (see also \cite{A2}).
\end{remark}

\subsection{Ramification data}
\label{sec:ram}

We will now consider the ramification data of the maps
\[
\phi_N:H^*\rightarrow X(N)
\]
and
\begin{equation}
\label{eqn:psiN} \psi_N: X(N)\rightarrow X(1)\cong \ppp^1.
\end{equation}
We will first study the characteristic zero case.  We recall some
facts from Shimura \cite{S}.

We call a matrix $g = \left(
\begin{array}{cc}
a & b\\
c & d
\end{array}
\right) \in SL(2,\zzz)$ {\it elliptic}, {\it parabolic}, or {\it
hyperbolic} if $|{\rm tr}(g)|$ is $<2$, $=2$ and $g\not= I$, or
$>2$, respectively. We call a point $z\in H^*$ {\it elliptic},
{\it parabolic}, or {\it hyperbolic} it it is fixed by such a
matrix. Parabolic points are also called {\it cuspidal points} or
{\it cusps}. Note that
 $g =
\left(
\begin{array}{cc}
a & b\\
c & d
\end{array}
\right)$, $c\not= 0$, fixes $z\in H$ if and only if

\begin{equation}
\label{eqn:fixes} z=\frac{a-d\pm \sqrt{(a+d)^2-4}}{2c}.
\end{equation}
(If $c=0$ then $g$ has no fixed points in $H$.)

If a group $G$ acts on a set $X$ then we denote by $G_x=\{g\in G\
|\ g(x)=x\}$ the stabilizer in $G$ of $x$.

\begin{lemma} (\cite{S})
\begin{itemize}
\item If $z\in H^*$ is elliptic then $z$ belongs to the
$SL_2(\zzz)$-orbit of $z_1=i$ or $z_2=\frac{1+\sqrt{3}i}{2}$.

\item The stabilizer in $SL_2(\zzz)$ of $z_1$ is

\[
SL_2(\zzz)_{z_1}= \langle \left(
\begin{array}{cc}
0 & 1\\
-1 & 0
\end{array}
\right) \rangle, \ \ \ \ \Gamma(N)_{z_1}=\{I\},
\]
and of $z_2$ is

\[
SL_2(\zzz)_{z_2}= \langle \left(
\begin{array}{cc}
0 & 1\\
-1 & 1
\end{array}
\right) \rangle, \ \ \ \ \Gamma(N)_{z_2}=\{I\}.
\]

\item If $z\in H^*$ is parabolic then $z$ belongs to the
$SL_2(\zzz)$-orbit of $\infty$. We have

\[
SL_2(\zzz)_{\infty}= \langle \left(
\begin{array}{cc}
1 & 1\\
0 & 1
\end{array}
\right) \rangle, \ \ \ \ \Gamma(N)_{\infty} =\langle \left(
\begin{array}{cc}
1 & N\\
0 & 1
\end{array}
\right) \rangle.
\]

\item There are no hyperbolic points in $H^*$,
\end{itemize}
\end{lemma}

The last part follows from (\ref{eqn:fixes}). The others are
proven in \cite{S}.

\begin{lemma}
(\cite{S}) The ramification index of $\psi_N$ at $P=\phi_N(z)$ is
$[\overline{\Gamma(1)}_z:\overline{\Gamma(N)}_z]$. The stabilizer
at $P$ is

\[
G_P=\overline{\Gamma(1)}_z/\overline{\Gamma(N)}_z.
\]
In particular, we have the following.

\begin{itemize}
\item If $z=z_1$ then

\[
G_P=H_1=\langle \left(
\begin{array}{cc}
0 & 1\\
-1 & 0
\end{array}
\right) \rangle
\]
is order $2$.

\item If $z=z_2$ then

\[
G_P=H_2=\langle \left(
\begin{array}{cc}
0 & 1\\
-1 & 1
\end{array}
\right) \rangle
\]
is order $3$.

\item If $z=\infty$ then

\[
G_P=H_3=\langle \left(
\begin{array}{cc}
1 & 1\\
0 & 1
\end{array}
\right) \rangle \cong \zzz/N\zzz
\]
is order $N$.
\end{itemize}

\end{lemma}

When $k=\ccc$, no other cyclic subgroups of $G$ occur as
stabilizers.

\begin{proposition}
\label{prop:spurious} Assume $k$ is an algebraically closed field of
good characteristic in the sense of Definition \ref{def:good}. If
$H$ is a cyclic subgroup of $G=Aut_k(X(N))=PSL_2(\zzz/N\zzz)$ which
is not of order $2$, $3$, or $N$ then $H$ is not the stabilizer of
any point on $X(N)$.  Furthermore, there is only one orbit of points
with stabilizer of order $2$, one orbit of points with stabilizer of
order $3$, and one orbit of points with stabilizer of order $N$.
\end{proposition}

\pf In the proof of Ritzenthaler \cite{R1} Proposition 1.3 it is
stated that $\phi_N$ is ramified over three points of indices $2$,
$3$ and $N$. The Hurwitz(-Riemann-Zeuthen) formula says that (see
for example Hartshorne, Corollary IV.2.4 \cite{H}):

\[
2\cdot g(X(N))-2=-2d+\deg(R),
\]
where $R$ is the ramification divisor

\[
R=\sum_{P\in X(N)} length(\Omega_{X(N)/X(1)})_P\cdot P.
\]
(Indeed, there are more explicit formulas for $\deg(R)$, since we
are in the tamely ramified case.) Since the genus of $X(N)$ and of
$X(1)$ do not depend on the characteristic, so we can conclude that
the ramification divisor of $X(N)/\ccc$ has the same degree as the
ramification divisor of $X(N)/k$. In other words, ``$\deg(R) \ {\rm
over}\ \ccc$'' equals ``$\deg(R) \ {\rm over}\ k$''. Since all the
summands in $\deg(R)$ are $\geq 0$, this is enough to show that no
other cyclic subgroups of $G$ arise as stabilizers over $k$. \qed

We call a cyclic subgroup of $G$ whose order is not $1$, $2$, $3$,
or $N$ {\bf spurious}. Spurious subgroups are not uncommon. Here is
a general construction. If $N$ is prime, so $\fff = \zzz/N\zzz$ is a
field, let $\eee/\fff$ denote a quadratic field with norm
$nm:\eee\rightarrow \fff$. The kernel of the norm map, $T=ker(nm)$
embeds into $SL_2(\fff)$. Moreover, since $\eee^\times$ is cyclic,
so is $T$, as well as all of its subgroups. It's known that $T$ has
order $N+1$ since the norm map is surjective in this case.

\begin{definition}
At each point $P \in X(N)$, the \textbf{ramification character}
$\theta_P$ is the character of the action of $G_P$ on the cotangent
space to $X(N)$ at $P$.
\end{definition}

\section{Representation theory of $PSL(2,N)$}
\label{sec:rep thy}

We first consider the representation theory of $G$ over $\ccc$,
following the treatment in \cite{FH}.

The group $PSL(2,N)$ has
$3+(N-1)/2$ conjugacy classes. Let $\varepsilon \in {\fff}$ be a
generator for the cyclic group ${\fff}^{\times}$.  Then each class
will have a representative of one of the following forms:

\[
\left(
\begin{array}{cc}
1 & 0\\
0 &1
\end{array}
\right),\
\left(
\begin{array}{cc}
x & 0\\
0 & x^{-1}
\end{array}
\right),\
\left(
\begin{array}{cc}
1 & 1\\
0 & 1
\end{array}
\right),\
\left(
\begin{array}{cc}
1 & \varepsilon\\
0 & 1
\end{array}
\right),\
\left(
\begin{array}{cc}
x & \varepsilon y\\
y & x
\end{array}
\right).
\]

The irreducible representations of $PSL(2,N)$ include the trivial
representation $\mathbf{1}$ and one irreducible $V$ of dimension
$N$. All but two of the others fall into two types:  representations
$W_{\alpha}$ of dimension $N+1$ (``principal series''), and
$X_{\beta}$ of dimension $N-1$ (``discrete series''). The first
type, $W_{\alpha}$, is labeled by a homomorphism $\alpha:
{\fff}^{\times} \to \ccc^{\times}$. The second type is indexed by a
homomorphism $\beta: T \to \ccc^{\times}$, where $T$ is a cyclic
subgroup of order $N+1$ of ${\fff}(\sqrt{\varepsilon})^{\times}$.
The characters of these are as follows:

{\footnotesize{
\[
\begin{array}[ht]{r||c|c|cc|c}
& \left( \begin{array}{cc} 1& 0 \\ 0 & 1 \end{array} \right) &
\left( \begin{array}{cc} x& 0 \\ 0 & x^{-1} \end{array} \right)&
\left( \begin{array}{cc} 1& 1 \\ 0 & 1 \end{array} \right)&
\left( \begin{array}{cc} 1& \varepsilon \\ 0 & 1 \end{array} \right)&
\left( \begin{array}{cc} x& \varepsilon y \\ y & x \end{array} \right)\\
\hline \hline
\mathbf{1} & 1 & 1 & 1 & 1 &1 \\
\hline X_{\beta} & N-1 & 0 & -1 & -1 &
-\beta(x+\sqrt{\varepsilon}y)-\beta(x-\sqrt{\varepsilon}y) \\
\hline V & N & 1 & 0 & 0 & -1 \\
\hline W_{\alpha} & N+1 & \alpha(x) + \alpha(x^{-1}) & 1 & 1 & 0 \\
\end{array}
\]
}}

Let $\tau$ denote a generator of $T$. Let $\zeta$ be a primitive
$N$th root of unity in $\ccc$.  Let $q$ and $q'$ be defined by

\begin{equation}
\label{eqn:qq'} q = \sum_{\left(\frac{a}{N}\right)=1} \zeta^a \mbox{
and } q' = \sum_{\left(\frac{a}{N}\right)=-1} \zeta^a,
\end{equation}
where the sums are over the quadratic residues and nonresidues
$\pmod N$, respectively. If $N \equiv 1$ mod 4, then the ``principal
series'' representation $W_{\alpha_0}$ corresponding to

\[
\begin{array}{ccc}
\alpha_0:{\fff}^{\times} & \to & \ccc^{\times} \\
\varepsilon & \mapsto & -1
\end{array}
\]
is not irreducible, but splits into two irreducibles $W'$ and $W''$,
each of dimension $(N+1)/2$.
Their characters satisfy:

\[
\begin{array}[ht]{r||c|c|cc|c}
& \left( \begin{array}{cc} 1& 0 \\ 0 & 1 \end{array} \right) &
\left( \begin{array}{cc} x& 0 \\ 0 & x^{-1} \end{array} \right)&
\left( \begin{array}{cc} 1& 1 \\ 0 & 1 \end{array} \right)& \left(
\begin{array}{cc} 1& \varepsilon \\ 0 & 1 \end{array} \right)&
\left( \begin{array}{cc} x& \varepsilon y \\ y & x \end{array} \right)\\
\hline \hline W' & \frac{N+1}{2} & \alpha_0(x) & 1+q
& 1+q' & 0 \\
\hline W'' & \frac{N+1}{2} & \alpha_0(x) & 1+q'
& 1+q & 0 \\
\end{array}
\]

Similarly, if $N \equiv 3$ mod 4, then the ``discrete series''
representation $X_{\beta}$ corresponding to

\[
\begin{array}{ccc}
\beta_0:T & \to & \ccc^{\times} \\
\tau & \mapsto & -1
\end{array}
\]
sending the generator $\tau$ to $-1$, splits into two
irreducibles $X'$ and $X''$, each
of dimension $(N-1)/2$. Their characters satisfy:

\[
\begin{array}[ht]{r||c|c|cc|c}
& \left( \begin{array}{cc} 1& 0 \\ 0 & 1 \end{array} \right) &
\left( \begin{array}{cc} x& 0 \\ 0 & x^{-1} \end{array} \right)&
\left( \begin{array}{cc} 1& 1 \\ 0 & 1 \end{array} \right)& \left(
\begin{array}{cc} 1& \varepsilon \\ 0 & 1 \end{array} \right)&
\left( \begin{array}{cc} x& \varepsilon y \\ y & x \end{array} \right)\\
\hline \hline X' & \frac{N-1}{2} & 0 & q
& q' & -\beta_0(x+y\sqrt{\varepsilon}) \\
\hline X'' & \frac{N-1}{2} & 0 & q'
& q & -\beta_0(x+y\sqrt{\varepsilon}) \\
\end{array}
\]
According to Janusz \cite{Ja}, the Schur index of each irreducible
representation of $G$ is $1$.

The action of the Galois group ${\cal G}$ on the irreducible
representations of $G$ is as follows.  The character values lie in
$\qqq(\mu)$, where $\mu$ is a primitive $m^{th}$ root of unity and
$m=N(N^2-1)/4$.  For each integer $j$ relatively prime to $m$, there
is an element $\sigma_j$ of the Galois group
${\cal{G}}=Gal(\qqq(\mu)/\qqq)$ taking $\mu$ to $\mu^j$.  This
Galois group element will act on representations by taking a
representation with character values $(a_1, \ldots, a_n)$ to a
representation with character values $(\sigma_j(a_1), \ldots,
\sigma_j(a_n))$. Representations with rational character values will
be fixed under this action.

Therefore the Galois group $\mathcal{G}$ will fix the trivial
representation and the $N$-dimensional representation $V$.  Its
action preserves the set of $N-1$-dimensional ``principal series''
representations $X_{\beta}$, and the set of $N+1$-dimensional
``discrete series'' representations $W_{\alpha}$. In the case $N
\equiv 1\pmod 4$, the Galois group will exchange the two
$(N+1)/2$-dimensional representations $W'$ and $W''$; if $N \equiv 3
\pmod 4$, the Galois group will exchange the two
$(N-1)/2$-dimensional representations $X'$ and $X''$.

\subsection{Brauer characters of some induced representations}

Let $k$ be a field of characteristic $\ell$, assume $k$ contains a 
primitive $|G|$-th root of unity\footnote{We really only need
for $k$ to contain the character values of the irreducible
$\overline{k}$-representations of $G$, where $\overline{k}$
denotes an algebraic closure of $k$.}, and assume that $\ell$
is good in the sense of (\ref{def:good}). Because $\ell$ does not
divide the order of $G$, the characters of $\ccc$-representations of
$G$ are the same as the Brauer characters of $k$-representations of
$G$. Some general remarks on the Brauer characters of induced
modular representations follow.

Let $G$ be a finite group, let $H$ be a subgroup.
Let $\sigma:H\rightarrow GL_m(k)$ be an $m$-dimensional
representation. Let

\[
\sigma^o(g)= \left\{
\begin{array}{ll}
\sigma(g),& g\in H,\\
0,&g\in G-H.
\end{array}
\right.
\]
If $G=g_1H\cup g_2H\cup ... \cup g_rH$ is a disjoint union into
cosets (where $g_i\in G$ and $g_1=1$) then let

\[
\pi = Ind_H^G\, \sigma (g)= (\sigma^o(g_i^{-1}gg_j))_{1\leq i,j\leq
r}.
\]
This $rm\times rm$ matrix is the induced representation of $\sigma$
to $G$. The trace of this representation is given by the $k$-valued
class function

\begin{equation}
\label{eqn:brauer} {\rm tr}\, \pi (g)= \sum_{x\in G/H}\ {\rm tr}\,
\sigma^o(x^{-1}gx),
\end{equation}
where we identify $G/H$ with $\{g_1,...,g_r\}$.

Fix an embedding $k^\times \hookrightarrow \ccc^\times$ (recall $k$
contains all the eigenvalues of each $\pi(g)$). Let $\chi_\sigma$
denote the Brauer character of $\sigma$ and let $\chi_\pi$ denote
the Brauer character of $\pi$.

\begin{proposition}
If $H$ has no elements of order $\ell={\rm char}(k)$ then

\[
\chi_\pi (g)= \sum_{x\in G/H}\ \chi_\sigma^o(x^{-1}gx),
\]
where $\chi_\sigma^o$ is the function $\chi_\sigma$ extended by $0$
off of $H$.
\end{proposition}

This is well-known and the proof of this is straightforward, using
only (\ref{eqn:brauer}) and the definitions, so omitted.

In the next few sections, we will use induced characters
to compute the $deg_{eq}(D)$ and $\tilde{\Gamma}_G$.
By the above proposition, such character computations 
will hold over $\ccc$ and in good positive characteristic. 
In particular, (\ref{eqn:JKrammod}) and all the examples
given in \S \ref{sec:examples} also
hold in this case.

\section{Induced characters}
\label{sec:inducedchars}

To compute the $PSL(2,N)$-module structure of the Riemann-Roch
space, we are interested in induced representations from the
non-spurious cyclic subgroups
\[
H_1=\langle \left(
\begin{array}{cc}
0 & 1\\
-1 & 0
\end{array}
\right) \rangle, \quad H_2=\langle \left(
\begin{array}{cc}
0 & 1\\
-1 & 1
\end{array}
\right) \rangle, \quad H_3=\langle \left(
\begin{array}{cc}
1 & 1\\
0 & 1
\end{array}
\right) \rangle
\]
of orders $2$, $3$, and $N$, respectively.  We compute these by
computing the restrictions of the irreducibles to these subgroups
and using Frobenius reciprocity.

\subsection{Induced characters from $H_1$}
\label{sec:IndH1}

To compute the restrictions, we must find the conjugacy classes of
$PSL(2,N)$ containing elements of these cyclic groups.  For $H_1$
and $H_2$, these depend on $N$.

Define the number $i$ as follows. When $N \equiv 1 \pmod 4$, let $i$
denote an element in ${\fff}^{\times}$ whose square is $-1$ (one can
take $i=\varepsilon^{(N-1)/4}$, where $\varepsilon$ is a generator
of $\fff^{\times}$). Then

\[
\left(
\begin{array}{cc}
0 & 1\\
-1 & 0
\end{array}
\right )\mbox{is conjugate to} \left(
\begin{array}{cc}
i & 0\\
0 & i^{-1}
\end{array}
\right ).
\]
When $N \equiv 3 \pmod 4$, there is no square root of -1 in $\fff$,
so we pass to the quadratic extension and let $i=x +
\sqrt{\varepsilon}y$ denote a square root of $-1$ in
${\fff}(\sqrt{\varepsilon})^{\times}$. (If $\tau$ is a generator of
the cyclic subgroup $T$ of order $N+1$ in
${\fff}(\sqrt{\varepsilon})^{\times}$, we can take
$i=\tau^{(N+1)/4}$.) Then

\[
\left(
\begin{array}{cc}
0 & 1\\
-1 & 0
\end{array}
\right )\mbox{is conjugate to} \left(
\begin{array}{cc}
x & \varepsilon y\\
y & x
\end{array}
\right ).
\]

Now let us compute the restrictions of the irreducible
representations of $PSL(2,N)$ to $H_1$.  We do this by examining the
relevant columns of the character table of $PSL(2,N)$ from section
\ref{sec:rep thy}, and comparing them to the character table of the
two element group $H_1$.

If $N \equiv 1 \pmod 4$, the relevant columns read

\[
\begin{array}[ht]{r||c|c}
& \left( \begin{array}{cc} 1& 0 \\ 0 & 1 \end{array} \right) &
\left( \begin{array}{cc} i& 0 \\ 0 & i^{-1} \end{array} \right)\\
\hline \hline
\mathbf{1} & 1 & 1  \\
\hline W' & \frac{N+1}{2} & \alpha_0(i) \\
\hline W'' & \frac{N+1}{2} & \alpha_0(i) \\
 \hline X_{\beta} & N-1 & 0  \\
\hline V & N & 1 \\
\hline W_{\alpha} & N+1 & \alpha(i) + \alpha(i^{-1}) \\
\end{array}
\]

Let $\theta_1$ be the nontrivial character of $H_1$.  From these
columns we see that the restricted representations are:
\begin{eqnarray*}
\Res_{H_1}^G \mathbf{1} & = & \mathbf{1}_{H_1} \\
\Res_{H_1}^{G} W' & = & \frac{N-1}{4} \mathbf{1}_{H_1} +
\frac{N-1}{4} \theta_1 +
\left\{ \begin{array}{c c l} \mathbf{1}_{H_1} & \mbox{if} & \alpha_0(i) = 1 \\
\theta_1 & \mbox{if} &  \alpha_0 (i) = -1 \end{array} \right. \\
\Res_{H_1}^{G} W'' & = & \mbox{ same } \\
\Res_{H_1}^{G} X_{\beta} & = & \frac{N-1}{2} \mathbf{1}_{H_1} + \frac{N-1}{2} \theta_1 \\
\Res_{H_1}^{G} V & = &  \frac{N+1}{2} \mathbf{1}_{H_1} + \frac{N-1}{2} \theta_1 \\
\Res_{H_1}^{G} W_{\alpha} & = & \frac{N-1}{2} \mathbf{1}_{H_1} +
\frac{N-1}{2} \theta_1 +
\left\{ \begin{array}{c c l} 2 \mathbf{1}_{H_1} & \mbox{if} & \alpha(i) = 1 \\
2 \theta_1 & \mbox{if} &  \alpha (i) = -1 \end{array} \right. \\
\end{eqnarray*}

The sign of $\alpha_0(i)$ depends on the equivalence class of $N
\pmod 8$. Since $\alpha_0(\varepsilon)=-1$ and $i =
\varepsilon^{(N-1)/4}$, we will have $\alpha_0(i)=1$ if $(N-1)/4$ is
even ($N \equiv 1 \pmod 8$) and $\alpha_0(i)=-1$ if $(N-1)/4$ is odd
($N \equiv 5 \pmod 8$).  For the induced representation of the
nontrivial character, we have, for $N \equiv 1 \pmod 8$,

{\footnotesize{
\begin{eqnarray*}
Ind_{H_1}^{G} \theta_1 & = & \frac{N-1}{4} (W' + W'') +
\frac{N-1}{2} \sum_{\beta} X_{\beta} + \frac{N-1}{2} V \\ & & +
\frac{N-1}{2} \sum_{\alpha(i)=1} W_{\alpha} + \frac{N+3}{2}
\sum_{\alpha(i)=-1} W_{\alpha},
\end{eqnarray*}
}} and for $N \equiv 5 \pmod 8$, {\footnotesize{
\begin{eqnarray*}
Ind_{H_1}^{G} \theta_1 & = & \frac{N+3}{4} (W' + W'') +
\frac{N-1}{2} \sum_{\beta} X_{\beta} + \frac{N-1}{2} V \\ && +
\frac{N-1}{2} \sum_{\alpha(i)=1} W_{\alpha} + \frac{N+3}{2}
\sum_{\alpha(i)=-1} W_{\alpha}.
\end{eqnarray*}
}}

Now let us compute the restricted representations for $N \equiv 3
\pmod 4$.  The relevant columns of the character table of $PSL(2,N)$
read

{\footnotesize{
\[
\begin{array}[ht]{r||c|c}
& \left( \begin{array}{cc} 1& 0 \\ 0 & 1 \end{array} \right) &
\left( \begin{array}{cc} x& \varepsilon y \\ y & x \end{array} \right)\\
\hline \hline
\mathbf{1} & 1 & 1  \\
\hline X' & \frac{N-1}{2} & -\beta_0(i) \\
\hline X'' & \frac{N-1}{2} & -\beta_0(i) \\
 \hline X_{\beta} & N-1 & -\beta(i) -\beta(i^{-1})  \\
\hline V & N & -1 \\
\hline W_{\alpha} & N+1 & 0 \\
\end{array}
\]
}} Therefore, the restricted representations are:

\begin{eqnarray*}
\Res_{H_1}^G \mathbf{1} & = & \mathbf{1}_{H_1} \\
\Res_{H_1}^{G} X' & = & \frac{N-3}{4} \mathbf{1}_{H_1} +
\frac{N-3}{4} \theta_1 +
\left\{ \begin{array}{c c l} \mathbf{1}_{H_1} & \mbox{if} & \beta_0(i) = -1 \\
\theta_1 & \mbox{if} &  \beta_0 (i) = 1 \end{array} \right. \\
\Res_{H_1}^{G} X'' & = & \mbox{ same } \\
\Res_{H_1}^{G} X_{\beta} & = & \frac{N-3}{2} \mathbf{1}_{H_1} +
\frac{N-3}{2} \theta_1 +
\left\{ \begin{array}{c c l} 2 \mathbf{1}_{H_1} & \mbox{if} & \beta(i) = -1 \\
2 \theta_1 & \mbox{if} &  \beta (i) = 1 \end{array} \right. \\
\Res_{H_1}^{G} V & = &  \frac{N-1}{2} \mathbf{1}_{H_1} + \frac{N+1}{2} \theta_1 \\
\Res_{H_1}^{G} W_{\alpha} & = & \frac{N+1}{2} \mathbf{1}_{H_1} +
\frac{N+1}{2} \theta_1 \\
\end{eqnarray*}

As before, the sign of $\beta_0(i)$ in the representations induced
from $X'$ and $X''$ depends on the equivalence class of $N \pmod 8$.
If $N \equiv 3 \pmod 8$, then $\beta_0(i)=-1$ and if $N \equiv 7
\pmod 8$, then $\beta_0(i) =1$.  We get the following induced
representations of the nontrivial character: if $N \equiv 3 \pmod
8$,

{\footnotesize{
\begin{eqnarray}
Ind_{H_1}^{G} \theta_1 & = & \frac{N-3}{4} (X' + X'') +
\frac{N+1}{2}
\sum_{\beta(i)=1} X_{\beta} + \frac{N-3}{2} \sum_{\beta(i)=-1} X_{\beta} \\
 & & + \frac{N+1}{2} V + \frac{N+1}{2} \sum_{\alpha} W_{\alpha}
\end{eqnarray}
}} and if $N \equiv 7 \pmod 8$,

{\footnotesize{
\begin{eqnarray}
Ind_{H_1}^{G} \theta_1 & = & \frac{N+1}{4} (X' + X'') +
\frac{N+1}{2}
\sum_{\beta(i)=1} X_{\beta} + \frac{N-3}{2} \sum_{\beta(i)=-1} X_{\beta} \\
 & & + \frac{N+1}{2} V + \frac{N+1}{2} \sum_{\alpha} W_{\alpha}.
\end{eqnarray} }}

\subsection{Induced characters from $H_2$}
\label{sec:IndH2}

Now we consider the induced characters from the cyclic group $H_2$
of order 3.  As in the case of $H_1$, the conjugacy classes of the
elements of $H_2$ depend on $N$.

We define a number $\omega$ as follows. When $N \equiv 1 \pmod 3$,
$\omega$ is a root of $x^2 + 1 = x$ in ${\fff}^{\times}$ (one can
take $\omega=\varepsilon^{(N-1)/6}$). Note

\[
\left(
\begin{array}{cc}
0 & 1\\
-1 & 1
\end{array}
\right )\mbox{is conjugate to} \left(
\begin{array}{cc}
\omega & 0\\
0 & \omega^{-1}
\end{array}
\right ),
\]
in $G$. When $N \equiv -1 \pmod 3$, $x + \sqrt{\varepsilon}y$ is a
root of $x^2 + 1 = x$ in ${\fff}(\sqrt{\varepsilon})^{\times}$ (one
can take $\omega=\tau^{(N+1)/6}$). Note

\[
\left(
\begin{array}{cc}
0 & 1\\
-1 & 0
\end{array}
\right )\mbox{is conjugate to} \left(
\begin{array}{cc}
x & \varepsilon y\\
y & x
\end{array}
\right ),
\]
in $G$. In any case, the nontrivial elements of $H_2$ are
$G$-conjugate.

Let $\theta_2$ and $\theta_2^2$ be the nontrivial characters of
$H_2$. Again, to compute $Ind_{H_2}^G \theta_2$ and $Ind_{H_2}^G
\theta_2^2$, we compute the restrictions of the irreducible
characters of $G$ and use Frobenius reciprocity.

If $N \equiv 1 \pmod 3$, the relevant columns of the character table
of $G=PSL(2,N)$ read

\[
\begin{array}[ht]{r||c|c}
& \left( \begin{array}{cc} 1& 0 \\ 0 & 1 \end{array} \right) &
\left( \begin{array}{cc} \omega & 0 \\ 0 & \omega^{-1} \end{array} \right)\\
\hline \hline
\mathbf{1} & 1 & 1  \\
\hline X' & \frac{N-1}{2} & 0 \\
\hline X'' & \frac{N-1}{2} & 0 \\
\hline W' & \frac{N+1}{2} & \alpha_0(\omega) \\
\hline W'' & \frac{N+1}{2} & \alpha_0(\omega) \\
 \hline X_{\beta} & N-1 & 0  \\
\hline V & N & 1 \\
\hline W_{\alpha} & N+1 & \alpha(\omega) + \alpha(\omega^{-1}) \\
\end{array}
\]
where the character table contains $W'$ and $W''$ if $N \equiv 1
\pmod 4$, and $X'$ and $X''$ if $N \equiv 3 \pmod 4$.  If $N \equiv
1 \pmod 4$, the value of $\alpha_0(\omega)$ is $(-1)^{(N-1)/6}$ and
$\frac{N-1}{6}$ is even, so $\alpha_0(\omega)=1$. Therefore, the
restricted representations are:

\begin{eqnarray*}
\Res_{H_2}^G \mathbf{1} & = & \mathbf{1}_{H_2} \\
\Res_{H_2}^{G} X' & = & \frac{N-1}{6} (\mathbf{1}_{H_2} + \theta_2 +
\theta_2^2) \\
\Res_{H_2}^{G} X'' & = & \mbox{ same } \\
\Res_{H_2}^G W' & = & \frac{N-1}{6} (\mathbf{1}_{H_2} + \theta_2 +
\theta_2^2) + \mathbf{1}_{H_2}\\
\Res_{H_2}^G W'' & = & \mbox{ same }\\
\Res_{H_2}^{G} X_{\beta} & = & \frac{N-1}{3} (\mathbf{1}_{H_2} +
\theta_2 + \theta_2^2) \\
\Res_{H_2}^{G} V & = &  \frac{N-1}{3} (\mathbf{1}_{H_2}
+ \theta_2 + \theta_2^2) + \mathbf{1}_{H_2}\\
\Res_{H_2}^{G} W_{\alpha} & = & \frac{N-1}{3} (\mathbf{1}_{H_2} +
\theta_2 + \theta_2^2) +
\left\{ \begin{array}{c c l} 2 \mathbf{1}_{H_2} & \mbox{if} & \alpha(\omega) = 1 \\
\theta_2 + \theta_2^2 & \mbox{if} &  \alpha (\omega) \neq 1  \end{array} \right. \\
Ind_{H_2}^{G} \theta_2 &=& Ind_{H_2}^{G} \theta_2^2 = \frac{N-1}{3}
\sum_{\beta} X_{\beta} + \frac{N-1}{3} V\\
&&\quad\quad\quad\quad\quad  + \frac{N-1}{3} \sum_{\alpha(\omega)=1}
W_{\alpha} + \frac{N+2}{3}
\sum_{\alpha(\omega) \neq 1} W_{\alpha} \\
&& \quad\quad\quad\quad\quad + \frac{N-1}{6} \left\{ \begin{array}{c
c l} W' + W'' & \mbox{if} & N \equiv 1 \pmod{12} \\ X' + X'' &
\mbox{if} & N \equiv 7 \pmod{12}
\end{array} \right. .
\end{eqnarray*}

Assume next $N \equiv 2 \pmod 3$. Again, the character table
contains $W'$ and $W''$ if $N \equiv 1 \pmod 4$ (so $N \equiv 5
\pmod{12}$), and $X'$ and $X''$ if $N \equiv 3 \pmod 4$ (so $N
\equiv 11 \pmod{12}$). The relevant columns of the character table
of $G$ are as follows. Note that for the irreducibles $X'$ and $X''$
(in the case $N \equiv 11 \pmod{12}$), the character value depends
on $\beta_0(\omega)$, where $\beta_0$ is the character of the cyclic
group of order $N+1$ sending a generator to $-1$.  In this case,
$\omega=\tau^{\frac{N+1}{6}}$, and $\frac{N+1}{6}$ is even, so
$\beta_0(\omega)=1$.

\[
\begin{array}[ht]{r||c|c}
& \left( \begin{array}{cc} 1& 0 \\ 0 & 1 \end{array} \right) &
\left( \begin{array}{cc} x & \varepsilon y \\ y & x \end{array} \right)\\
\hline \hline
\mathbf{1} & 1 & 1  \\
\hline X' & \frac{N-1}{2} & -1 \\
\hline X'' & \frac{N-1}{2} & -1 \\
\hline W' & \frac{N+1}{2} & 0 \\
\hline W'' & \frac{N+1}{2} & 0 \\
 \hline X_{\beta} & N-1 & -\beta(\omega) - \beta(\omega^{-1}) \\
\hline V & N & -1 \\
\hline W_{\alpha} & N+1 &  0\\
\end{array}
\]
Therefore, we see that when $N \equiv 2 \pmod 3$, the restricted
representations are:
\begin{eqnarray*}
\Res_{H_2}^G \mathbf{1} & = & \mathbf{1}_{H_2} \\
\Res_{H_2}^{G} X' & = & \frac{N-5}{6} (\mathbf{1}_{H_2} + \theta_2 +
\theta_2^2) + \theta_2 + \theta_2^2 \\
\Res_{H_2}^{G} X'' & = & \mbox{ same } \\
\Res_{H_2}^G W' & = & \frac{N+1}{6} (\mathbf{1}_{H_2} + \theta_2 +
\theta_2^2)\\
\Res_{H_2}^G W'' & = & \mbox{ same }\\
\Res_{H_2}^{G} X_{\beta} & = & \frac{N-5}{3} (\mathbf{1}_{H_2} +
\theta_2 + \theta_2^2) + \left\{ \begin{array}{c c l} 2( \theta_2 + \theta_2^2) & \mbox{if} & \beta(\omega) =1 \\
 2 \mathbf{1}_{H_2} + \theta_2 + \theta_2^2 & \mbox{if} & \beta(\omega) \neq 1 \end{array} \right. \\
\Res_{H_2}^{G} V & = &  \frac{N-2}{3} (\mathbf{1}_{H_2} + \theta_2 + \theta_2^2) + \theta_2 + \theta_2^2 \\
\Res_{H_2}^{G} W_{\alpha} & = & \frac{N+1}{3} (\mathbf{1}_{H_2} +
\theta_2 + \theta_2^2) \\
\end{eqnarray*}

And the induced representations for $N \equiv 2 \pmod 3$ are:
\begin{eqnarray*}
Ind_{H_2}^{G} \theta_2 = Ind_{H_2}^{G} \theta_2^2 &=& \frac{N+1}{3}
\sum_{\beta} X_{\beta} + \frac{N-2}{3} \sum_{\beta} X_{\beta} +
\frac{N+1}{3}\sum_{\alpha}
W_{\alpha} \\
&&+\frac{N+1}{3} V + \frac{N+1}{6} \left\{ \begin{array}{c c l} W' +
W'' & \mbox{if} &
N \equiv 5 \pmod{12} \\
X' + X'' & \mbox{if} & N \equiv 11 \pmod{12} \end{array} \right.
\end{eqnarray*}

\subsection{Induced characters from $H_3$}
\label{sec:IndH3}

Now we consider the induced characters from $H_3$, the cyclic group
of order $N$.  For the restrictions of the irreducible characters of
$PSL(2,N)$ to $H_3$, the relevant columns of the character table of
$G=PSL(2,N)$ are as follows.  As before, the representations $W'$
and $W''$ appear if and only if
$N \equiv 1 \pmod 4$, and $X'$ and $X''$ appear
if and only if $N \equiv 3 \pmod 4$.

\[
\begin{array}[ht]{r||c|c c}
& \left( \begin{array}{cc} 1& 0 \\ 0 & 1 \end{array} \right) &
\left( \begin{array}{cc} 1& 1 \\ 0 & 1 \end{array} \right)&
\left( \begin{array}{cc} 1& \varepsilon \\ 0 & 1 \end{array} \right)\\
\hline \hline
\mathbf{1} & 1 & 1 & 1 \\
\hline X' & \frac{N-1}{2} & q & q' \\
\hline X'' & \frac{N-1}{2} & q' & q \\
\hline W' & \frac{N+1}{2} & 1+q & 1+q' \\
\hline W'' & \frac{N+1}{2} & 1+q'& 1+q \\
 \hline X_{\beta} & N-1 & -1 & -1  \\
\hline V & N & 0 & 0 \\
\hline W_{\alpha} & N+1 & 1 & 1 \\
\end{array}
\]
Recall
\[
q = \sum_{\left(\frac{a}{N}\right)=1} \zeta^a \mbox{ and } q' =
\sum_{\left(\frac{a}{N}\right)=-1} \zeta^a.
\]

Now let $\theta_3$ be the nontrivial character of $H_3$ such that
\[
\theta_3 \left( \begin{array}{cc} 1 & 1 \\ 0 & 1 \end{array} \right)
= \zeta.
\]
The restricted representations are:
\begin{eqnarray*}
\Res_{H_3}^G \mathbf{1} & = & \mathbf{1}_{H_3} \\
\Res_{H_3}^{G} X' & = & \sum_{\left(\frac{k}{n}\right)=1} \theta_3^k \\
\Res_{H_3}^{G} X'' & = & \sum_{\left(\frac{k}{n}\right)=-1} \theta_3^k \\
\Res_{H_3}^G W' & = & \mathbf{1}_{H_3} + \sum_{\left(\frac{k}{n}\right)=1} \theta_3^k\\
\Res_{H_3}^G W'' & = & \mathbf{1}_{H_3} + \sum_{\left(\frac{k}{n}\right)=-1} \theta_3^k\\
\Res_{H_3}^{G} X_{\beta} & = & \sum_{k=1}^{N-1} \theta_3^k \\
\Res_{H_3}^{G} V & = & \mathbf{1}_{H_3} + \sum_{k=1}^{N-1} \theta_3^k \\
\Res_{H_3}^{G} W_{\alpha} & = & 2\mathbf{1}_{H_3} + \sum_{k=1}^{N-1} \theta_3^k \\
\end{eqnarray*}
Threrefore,

{\footnotesize{
\[
Ind_{H_3}^{G} \theta_3^k = \sum_{\beta} X_{\beta} + V +
\sum_{\alpha}
W_{\alpha} \\
+ \left\{ \begin{array}{ll}
W', &N \equiv 1 \pmod 4, \ k\ {\rm a\ quad.\ residue\ mod\ }N,\\
X', &N \equiv 3 \pmod 4, \ k\ {\rm a\ quad.\ residue\ mod\ }N,\\
W'', &N \equiv 1 \pmod 4, \ k\ {\rm not\ a\ quad.\ residue\ mod\ }N,\\
X'', &N \equiv 3 \pmod 4, \ k\ {\rm not\ a\ quad.\ residue\ mod\ }N.
\end{array} \right.
\]
}}

\section{Ramification module}
\label{sec:rammod}

For the remainder of this paper, let $G^*$ denote the set of
equivalence classes of irreducible representations of $G$.  We often
abuse notation by using the same symbol for an element of $G^*$, a
representation in its equivalence class, and the character of such a
representation. Over a field $k$ of good positive characteristic, we
abuse notation by using $G^*$ instead to represent Brauer characters
of $G$ over $k$.

Define the module $\Gamma_G$ by

\[
\Gamma_G=\sum_{P\in X(N)} Ind_{G_P}^G\left( \sum_{\ell =1}^{e_P-1}
\ell\theta_P^\ell \right),
\]
where $\theta_P$ is the ramification character at a point $P$.
Thanks to Nakajima \cite{N}, it is known that there is a $G$-module
$\tilde{\Gamma}_G$ such that $[\Gamma_G]=|G|\cdot[
\tilde{\Gamma}_G]$ (see also \cite{B}). We call $\tilde{\Gamma}_G$
the {\it ramification module}.

In the case of the modular curve $X(N)$, let $\pi_{\theta_P}$ denote
the induced representation from a character $\theta_P$.  Then by
definition,

\begin{equation}
\tilde{\Gamma}_G= \frac{\overline{R}_2[\pi_{\theta_1}]
+\overline{R}_3 ([\pi_{\theta_2}]+2[\pi_{\theta_2^2}])
+\overline{R}_N ([\pi_{\theta_3}] +...+
(N-1)[\pi_{\theta_3^{N-1}}])}{|G|},
\end{equation}
where $\overline{R}_i$ is the number of $P\in X(N)$ with $G_P=H_i$.
Recall from section \ref{sec:ram} that

\[
\overline{R}_2=\frac{|G|}{2},\ \ \ \overline{R}_3=\frac{|G|}{3},\ \
\ \overline{R}_N=\frac{|G|}{N},
\]
so that in fact
\begin{equation}
\label{eqn:ram_mod} \tilde{\Gamma}_G= \frac{1}{2}[\pi_{\theta_1}] +
\frac{1}{3} ([\pi_{\theta_2}]+2[\pi_{\theta_2^2}]) +\frac {1}{N}
([\pi_{\theta_3}] +...+ (N-1)[\pi_{\theta_3^{N-1}}]).
\end{equation}

This can easily be computed using the induced characters from
section \ref{sec:inducedchars}. To combine the results of sections
\ref{sec:IndH1} and \ref{sec:IndH2}, we must look at the congruence
class of $N$ modulo $24$. We have

{\small{ $N\equiv 1\pmod {24}\implies N\equiv 1\pmod {4}, N\equiv
1\pmod {8}, N\equiv 1\pmod {12}$,

$N\equiv 5\pmod {24}\implies N\equiv 1\pmod {4}, N\equiv 5\pmod {8},
N\equiv 5\pmod {12}$,

$N\equiv 7\pmod {24}\implies N\equiv 3\pmod {4}, N\equiv 7\pmod {8},
N\equiv 7\pmod {12}$,

$N\equiv 11\pmod {24}\implies N\equiv 3\pmod {4}, N\equiv 3\pmod
{8}, N\equiv 11\pmod {12}$,

$N\equiv 13\pmod {24}\implies N\equiv 1\pmod {4}, N\equiv 5\pmod
{8}, N\equiv 1\pmod {12}$,

$N\equiv 17\pmod {24}\implies N\equiv 1\pmod {4}, N\equiv 1\pmod
{8}, N\equiv 5\pmod {12}$,

$N\equiv 19\pmod {24}\implies N\equiv 3\pmod {4}, N\equiv 3\pmod
{8}, N\equiv 7\pmod {12}$,

$N\equiv 23\pmod {24}\implies N\equiv 3\pmod {4}, N\equiv 7\pmod
{8}, N\equiv 11\pmod {12}$. }} Let ${\cal Q}$ denote the set of
quadratic residues $\pmod N$, let ${\cal N}$ denote the set of
quadratic non-residues $\pmod N$ and let

\[
S_{\cal Q} =\sum_{\ell\in {\cal Q}}\ell,\ \ \ \ S_{\cal
N}=\sum_{\ell\in {\cal N}}\ell = \frac{N(N-1)}{2}-S_{\cal Q}.
\]
If $N \equiv 1 \pmod 4$, then

\begin{equation}
\label{eqn:S_q}
S_{\cal Q}=S_{\cal N}=\frac{N(N-1)}{4}.
\end{equation}
If $N \equiv 3 \pmod 4$, then $S_{\cal Q}$ and $S_{\cal N}$ are not
equal.

We prove the following result.  An empty sum is by convention zero.

\begin{theorem}
\label{thrm:main} We have the following decomposition of the
ramification module:

\begin{itemize}
\item
If $N\equiv 1\pmod {24}$, let $m=\frac{13N-13}{12}$.  Then
\begin{eqnarray*}
\tilde{\Gamma}_G & = &\frac{m}{2}W'+\frac{m}{2}W''+mV +
m\sum_{\beta} X_{\beta} \\
&&+ m\sum_{\alpha(i)=\alpha(\omega)=1} W_{\alpha}
+(m+1)\sum_{\alpha(i)=1,\alpha(\omega)\not= 1} W_{\alpha}\\
&&+(m+1)\sum_{\alpha(i)\not= 1,\alpha(\omega) = 1} W_{\alpha}
+(m+2)\sum_{\alpha(i)\not= 1,\alpha(\omega)\not= 1} W_{\alpha}.
\end{eqnarray*}

\item
If $N\equiv 5\pmod {24}$, let $m=\frac{13N-5}{12}$. Then
\begin{eqnarray*}
\tilde{\Gamma}_G & = & \frac{m+1}{2}W' + \frac{m+1}{2}W'' + mV \\
&&+ m\sum_{\beta(\omega)=1} X_{\beta}
+(m-1)\sum_{\beta(\omega)\not=1} X_{\beta}\\
&& + m\sum_{\alpha(i)=1} W_{\alpha} +(m+1)\sum_{\alpha(i)\not= 1}
W_{\alpha}.
\end{eqnarray*}

\item
If $N\equiv 7\pmod {24}$, let $m=\frac{13N-7}{12}$. Then
\begin{eqnarray*}
\tilde{\Gamma}_G & = & +(\frac{7N-1}{24}+\frac{S_{\cal Q}}{N}))X'+
(\frac{7N-1}{24}+\frac{S_{\cal N}}{N}))X'' + mV\\
&&+m\sum_{\beta(i)=1} X_{\beta} +(m-1)\sum_{\beta(i)=-1}X_{\beta}\\
&&+ m\sum_{\alpha(\omega)=1} W_{\alpha}
+(m+1)\sum_{\alpha(\omega)=-1} W_{\alpha}.
\end{eqnarray*}

\item
If $N\equiv 11\pmod {24}$, let $m=\frac{13N+1}{12}$. Then
\begin{eqnarray*}
\tilde{\Gamma}_G & = & (\frac{7N-5}{24}+\frac{S_{\cal Q}}{N})X'
+(\frac{7N-5}{24}+\frac{S_{\cal N}}{N})X'' +mV +
m\sum_{\alpha} W_{\alpha}\\
&&+m\sum_{\beta(i)=\beta(\omega)=1} X_{\beta}
+(m-1)\sum_{\beta(i)=1,\beta(\omega)\not= 1} X_{\beta}\\
&&+(m-1)\sum_{\beta(i)=-1,\beta(\omega)=1} X_{\beta}
+(m-2)\sum_{\beta(i)=-1,\beta(\omega)\not= 1} X_{\beta}.
\end{eqnarray*}

\item
If $N\equiv 13\pmod {24}$, let $m=\frac{13N-13}{12}$. Then
\begin{eqnarray*}
\tilde{\Gamma}_G & = & \frac{m}{2}W' +\frac{m}{2}W'+ mV
+m\sum_{\beta}X_{\beta}\\
&&+  m\sum_{\alpha(i)=\alpha(\omega)=1} W_{\alpha}
+(m+1)\sum_{\alpha(i)=1,\alpha(\omega)\not= 1} W_{\alpha}\\
&&+(m+1)\sum_{\alpha(i)\not= 1,\alpha(\omega) = 1} W_{\alpha}
+(m+2)\sum_{\alpha(i)\not= 1,\alpha(\omega)\not= 1} W_{\alpha}.
\end{eqnarray*}

\item
If $N\equiv 17\pmod {24}$, let $m=\frac{13N-5}{12}$.  Then
\begin{eqnarray*}
\tilde{\Gamma}_G & = & \frac{m}{2}W'+\frac{m}{2}W'' +mV\\
&&+m\sum_{\beta(\omega)=1} X_{\beta}+ (m-1)\sum_{\beta(\omega)\not=
1} X_{\beta} \\
&&+m\sum_{\alpha(i)=1} W_{\alpha} +(m+1)\sum_{\alpha(i)\not= 1}
W_{\alpha}.
\end{eqnarray*}

\item
If $N\equiv 19\pmod {24}$, let $m=\frac{13N-7}{12}$. Then
\begin{eqnarray*}
\tilde{\Gamma}_G & = & (\frac{7N-13}{24}+\frac{S_{\cal Q}}{N})X' +
(\frac{7N-13}{24}+\frac{S_{\cal N}}{N})X'' + mV\\
&&+ m\sum_{\beta(i)=1} X_{\beta}+(m-1)\sum_{\beta(i)=-1} X_{\beta}\\
&&+m\sum_{\alpha(\omega)=1} W_{\alpha} +(m+1)\sum_{\alpha(\omega)=1}
W_{\alpha}.
\end{eqnarray*}

\item
If $N\equiv 23\pmod {24}$, let $m=\frac{13N+1}{12}$. Then
\begin{eqnarray*}
\tilde{\Gamma}_G & = & (\frac{7N+7}{24}+{\frac{S_{\cal{Q}}}{N}})X'
+(\frac{7N+7}{24}+\frac{S_{\cal N}}{N})X''+mV +m\sum_{\alpha} W_{\alpha}\\
&&+m\sum_{\beta(i)=\beta(\omega)=1} X_{\beta}
+(m-1)\sum_{\beta(i)=1,\beta(\omega)=\not= 1} X_{\beta}\\
&&+(m-1)\sum_{\beta(i)=-1,\beta(\omega)=1} X_{\beta}
+(m-2)\sum_{\beta(i)=-1,\beta(\omega)\not= 1} X_{\beta}.
\end{eqnarray*}

\end{itemize}
\end{theorem}

\pf It is straightforward to verify the theorem case-by-case using
(\ref{eqn:ram_mod}).  As an intermediate step, we tabulate the
multiplicities of the irreducible representations in the induced
representations which comprise $\Gamma_G$. These tables are included
in Appendix \S \ref{sec:appendix}.  We sketch the case $N\equiv
1\pmod{24}$.

We begin with the multiplicity of $V$.  According to
(\ref{eqn:ram_mod}) and the $N\equiv 1\pmod{24}$ table in Appendix
\S \ref{sec:appendix}, the multiplicity of $V$ in
$\tilde{\Gamma}_G$. is
\[
\frac{1}{2} \frac{N-1}{2} + \frac{1}{3} N-1 +
\frac{1}{N}\frac{N(N-1)}{2}= \frac{13N-13}{12} = m.
\]
The multiplicity of an $X_\beta$ (or a $W_\alpha$ with
$\alpha(i)=\alpha(\omega)=1$) in $\tilde{\Gamma}_G$ is the same,
since they have the same rows in the table as $V$.  Next consider
the multiplicities of $W'$ and $W''$.  Since $S_{\cal Q}=S_{\cal
N}=\frac{N(N-1)}{4}$, these rows in the table are exactly half of
the $V$ row, so the multiplicities of $W'$ and $W''$ are both
\[
\frac{1}{2} \frac{N-1}{4} + \frac{1}{3} \frac{N-1}{2} +
\frac{1}{N}\frac{N(N-1)}{4}= \frac{13N-13}{24} = \frac{m}{2}.
\]

Next consider a $W_\alpha$ with $\alpha(i)=1$, $\alpha(\omega)\not=
1$. This row of the table differs from the $V$ row in the second
column, where it is $N+2$ instead of $N-1$.  Since this entry is
divided by $3$ in the computation, the multiplicity in
$\tilde{\Gamma}_G$ of such a $W_\alpha$ will be one more than the
multiplicity of $V$, or $m+1$.  More directly,
\[
\frac{1}{2} \frac{N-1}{2} + \frac{1}{3} N+2 +
\frac{1}{N}\frac{N(N-1)}{2}= \frac{13N-13}{12} = m+1.
\]

Similarly, consider a $W_\alpha$ with $\alpha(i)\not= 1$,
$\alpha(\omega) = 1$.  This row of the table differs from the $V$
row in the first column, where again the net result is to add one to
the multiplicity.  In this case,
\[
\frac{1}{2} \frac{N+3}{2} + \frac{1}{3} N-1 +
\frac{1}{N}\frac{N(N-1)}{2}= \frac{13N-13}{12} =m+1.
\]
as desired. Finally, consider a $W_\alpha$ with $\alpha(i)\not= 1$,
$\alpha(\omega) \not= 1$. We see that both of the first two columns
of the table differ in this row from the $V$ row, for a net addition
of $2$ to the multiplicity.  In this case,
\[
\frac{1}{2} \frac{N+3}{2} + \frac{1}{3} N+2 +
\frac{1}{N}\frac{N(N-1)}{2}= \frac{13N-13}{12} = m+2.
\]
This completes the verification of the theorem in the case $N\equiv
1\pmod{24}$.  A similar analysis can be used on the other cases.\qed

If $\tilde{\Gamma}_G$ has a $\qqq[G]$-module structure, it may be
computed more simply \cite{JK}.   In this case, the formula for it
is

\begin{equation}
\label{eqn:JKrammod} \tilde{\Gamma}_G=\bigoplus_{\pi\in G^*} \left[
\sum_{\ell=1}^L ({\rm dim}\, \pi -{\rm dim}\, (\pi^{H_\ell}))
\frac{R_\ell}{2}\right]\pi
\end{equation}
where  $\{H_1,...,H_L\}$ represent the set of conjugacy classes of
cyclic subgroups of $G$). If $\tilde{\Gamma}_G$ is not ${\cal
  G}=Gal(\qqq(G)/\qqq)$-invariant,
so there is no $\qqq[G]$-module structure, then

\begin{equation}
\label{eqn:JKGalclosure}
\bigoplus_{\pi\in G^*} |{\cal G}| \cdot
\left[ \sum_{\ell=1}^L ({\rm dim}\, \pi -{\rm dim}\, (\pi^{H_\ell}))
\frac{R_\ell}{2}\right]\pi
\end{equation}
yields the $G$-module decomposition of the Galois-closure of
$\tilde{\Gamma}_G$ (i.e., the smallest  ${\cal G}$-invariant
$G$-module containing $\tilde{\Gamma}_G$). See the remark after
Corollary 6 in \cite{JK}. Also, compare this with Lemma 3 in
\cite{KP} (we thank Dipendra Prasad for providing us with this
reference).

This motivates the following, as stated in the introduction.

\begin{theorem}
\label{thrm:2}
For $N>5$ prime, the ramification module of
$X(N)$ over $X(1)$ is ${\cal G}$-invariant if and only if $N\equiv 1
\pmod 4$.
\end{theorem}

\pf.
As noted in \S \ref{sec:rep thy}, the action of ${\cal G}$ is as
follows:

\begin{itemize}

\item
{\it On $X'$ and $X''$ or $W'$ and $W''$}:  The induced characters
from $H_1$ and $H_2$ are invariant under this action. The induced
characters from $H_3$ are not; the action exchanges $\theta^k$ where
$k$ is a quadratic residue $\pmod N$ with $\theta^{k'}$ where $k'$
is not a quadratic residue $\pmod N$.  However the sum appearing in
the ramification module
\[
Ind_{G_P}^G\left( \sum_{\ell =1}^{e_P-1} \ell\theta_P^\ell \right)
\\
= \sum_{\ell=1}^{N-1} \ell Ind_{H_3}^G \theta^{\ell}
\]
will include
\[
\sum_{\left(\frac{\ell}{N} \right)=1} \ell W' +
\sum_{\left(\frac{\ell}{N} \right)=-1} \ell W''
\]
if $N \equiv 1 \pmod 4$, or
\[
\sum_{\left(\frac{\ell}{N} \right)=1} \ell X' +
\sum_{\left(\frac{\ell}{N} \right)=-1} \ell X''
\]
if $N \equiv 3 \pmod 4$.
Here $\left(\frac{*}{N} \right)$ denotes the Legendre symbol.
If $N \equiv 1 \pmod 4$, then the sum of
the quadratic residues is the same as the sum of the nonresidues, so
the multiplicities of $W'$ and $W''$ are the same. If $N \equiv 3
\pmod 4$, then the sum of the quadratic residues is not the same as
the sum of the quadratic nonresidues, so these multiplicities are
different.  Therefore this part of the ramification module module is
Galois invariant if and only if $N \equiv 1 \pmod 4$.

\item
{\it On $X_\beta$}: The multiplicity of a discrete series
representation $X_{\beta}$ depends only on $\beta(i)$ and
$\beta(\omega)$. Observe that $\mathcal{G}$ fixes elements of
$\qqq$, therefore cannot exchange $X_{\beta}$ having $\beta(i)=1$
with an $X_{\beta}$ having $\beta(i)=-1$, or $X_{\beta}$ having
$\beta(\omega)=1$ with $X_{\beta}$ having $\beta(\omega) \neq 1$.
Therefore this part of the ramification module is Galois invariant.

\item
{\it On $W_\alpha$}:  Similarly, the multiplicity of a principal
series representation $W_{\alpha}$ depends only on $\alpha(i)$ and
$\alpha(\omega)$.  Again, $\mathcal{G}$ will not exchange these, so
this part of the module is Galois invariant.
\end{itemize}

Thus we see that $\tilde{\Gamma}_G$ is Galois-invariant if and only
if $N \equiv 1 \pmod 4$. \qed

In this case, we can use formula (\ref{eqn:JKrammod}) to compute the
ramification module directly from the restricted representations in
section \ref{sec:inducedchars}, and get the same result with fewer
steps. If $N \equiv 3 \pmod 4$, formula (\ref{eqn:JKGalclosure})
will yield the Galois closure of the ramification module.

\section{Equivariant degree and Riemann-Roch space}
\label{sec:equivdeg}

Now we will define and compute the equivariant degree of a
$G$-invariant divisor.  This, together with Borne's formula
(\ref{eqn:Borne}) will allow us to compute the $G$-module structure
of the Riemann-Roch space $L(D)$.

Fix a point $P\in X(N)$ and let $D$ be a divisor on $X(N)$ of the
form

\[
D=\frac{r}{e_P}\sum_{g\in G}g(P)=r\sum_{g\in G/G_P}g(P),
\]
where $G_P$ denotes the stabilizer in $G$ of $P$ and $e_P=|G_P|$
denotes the ramification index at $P$.

When $r=1$ such a divisor is called a {\it reduced orbit} by
Borne \cite{B}.
When $r$ is a multiple of $e_P$ then $D$ is the pull-back of a divisor
on $X(1)$ via $\psi_N$ in (\ref{eqn:psiN}).

When $D$ is as above,
the {\it equivariant degree} of $D$ is the virtual representation

\[
deg_{eq}(D)=Ind_{G_P}^G\ \left(\sum_{\ell =1}^r
\theta_P^{-\ell}\right),
\]
where $\theta_P$ is the ramification character of $X(N)$ at $P$ (a
character of $G_P$). In general, the equivariant degree is additive
on disjointly supported divisors. See for example Borne \cite{B} for
more details.

\begin{example}
Consider a hyperelliptic curve $X$ for which every ramification
point $P$ has $e_P=2$. Let $D$ denote the sum of all the
ramification points of $X$. Then (by definition)
$deg_{eq}(D)=\Gamma_G=2 \tilde{\Gamma}_G$.
\end{example}

On the modular curve $X(N)$, the results of section \ref{sec:ram}
tell us that there are only four types of reduced orbits to
consider:  the stabilizer $G_P$ of a point $P$ in the support of $D$
may have order $1$, $2$, $3$, or $N$.  Let $D_0$, $D_1$, $D_2$, and
$D_3$ denote reduced orbits of each type.  We compute the
equivariant degree of $rD_i$ in each case.

\begin{description}

\item[Case 1]: $G_P$ is trivial.
If $G_P$ has order $1$, then the ramification character is trivial.
$D_0$ is an entire orbit, in fact the pullback of a divisor on
$X(1)$.   The equivariant degree of $rD_0$ is $r \cdot k[G]$, a
multiple of the regular representation $k[G]$ of $G$.

\item[Case 2]: $G_P \cong H_1$.
If $G_P$ has order $2$, then the equivariant degree of $D_1$ is
$Ind_{G_P}^G \theta_1$ (which was computed in section
\ref{sec:IndH1}).  Even multiples of $D_1$ will be pullbacks from
$X(1)$, so $\deg_{eq}(2rD_1) = r \cdot k[G]$, and
$\deg_{eq}((2r+1)D_1) = r \cdot k[G] + Ind_{G_P}^G \theta_1$.

\item[Case 3]: $G_P \cong H_2$.  The
equivariant degree of $D_2$ is $Ind_{G_P}^G \theta_2^2$, which was
computed in section \ref{sec:IndH2}.  Note that $Ind_{G_P}^G
\theta_2^2= Ind_{G_P}^G \theta_2$, so $\deg_{eq}(2D_2)= 2
Ind_{G_P}^G \theta_2$.  But $3 D_2$ is a pullback from $X(1)$, so
$\deg_{eq}(3D_2)= k[G]$, and more generally $\deg_{eq}(3rD_2)= r
\cdot k[G]$, $\deg_{eq}((3r+1)D_2)= r \cdot k[G] + Ind_{G_P}^G
\theta_2$, and $\deg_{eq}((3r+2)D_2)= r \cdot k[G] + 2 Ind_{G_P}^G
\theta_2$.

\item[Case 4]:  $G_P \cong H_3$.  This case is more complicated
because $Ind_{H_3}^{G} \theta_3^k$ depends on whether $k$ is a
quadratic residue modulo $N$ (see section \ref{sec:IndH3}).  If $1
\leq r \leq N-1$, let $q_{r}$ denote the number of quadratic
residues modulo $N$ in the set $N-r, \ldots, N-1$, and let $n_{r}$
denote the number of quadratic nonresidues modulo $N$ in the same
set.  Then the equivariant degree $\deg_{eq}(r D_3)$ will be

{\footnotesize{
\[
\deg_{eq}(r D_3) = r \left( \sum_{\beta} X_{\beta} + V +
\sum_{\alpha}
W_{\alpha} \right) \\
+ \left\{ \begin{array}{ll}
q_r W' + n_r W'', &N \equiv 1 \pmod 4\\
q_r X' + n_r X'', &N \equiv 3 \pmod 4.
\end{array} \right.
\]
}}

\end{description}

Now we would like to compute the $G$-module structure of the
Riemann-Roch space $L(D)$ for a non-special $G$-invariant divisor
$D$. First, let us consider which $G$-invariant divisors are
non-special.  To be non-special, it is sufficient to have $deg D >
2g-2$, where
\[
g=1+\frac{(N-6)(N^2-1)}{24}
\]
is the genus of $X(N)$, so
\[
2g-2 = \frac{(N-6)(N^2-1)}{12} = \frac{(N-6)}{6N} |G|.
\]
The reduced orbits $D_0$, $D_1$ and $D_2$ have degrees $|G|$,
$|G|/2$ and $|G|/3$, respectively, which are greater than $2g-2$ for
any $N$, so these are always non-special.  However, the reduced
orbit $D_3$ has degree $|G|/N$, which is greater than $2g-2$ only if
$N<12$.  A multiple of this reduced orbit, $rD_3$ will be
non-special if $r > (N-6)/6$.

Once the equivariant degree of a non-special divisor has been
computed, it is simple to combine this with the ramification module,
using Borne's formula, to find the $G$-module structure of the
Riemann-Roch space $L(D)$. The number of cases becomes quite
cumbersome, but in each case there will be a number $M$ such that
all but the two smallest nontrivial irreducibles (either $W'$ and
$W''$ or $X'$ and $X''$) will have multiplicity either $M$, $M+1$,
$M+2$, or $M+3$, and the multiplicities of the two smallest
irreducibles will be half of one of these numbers.

There is one special case worth mentioning:  if $D$ is a pull-back
of an effective divisor $\overline{D}$ on $X(1)$ via $\psi_N$, then
it will be non-special; moreover in this case
\[
[\deg_{eq}(D)]=\deg(D_0)[k[G]],
\]
and we will have the following decomposition of the Riemann-Roch
space
\[
[L(D)]=(1+\deg(\overline{D}))[k[G]]-[\tilde{\Gamma}_G],
\]
where $\tilde{\Gamma}_G$ is as above.

\section{Examples}
\label{sec:examples}

The computer algebra system \cite{GAP} computes information about
$PSL(2,N)$; one can use it to compute character tables, induced
characters and Schur inner products.  In the examples of $X(7)$ and
$X(11)$ below, we use GAP to explicitly compute the $G$-module
structure of the ramification module and some Riemann-Roch spaces,
in the cases $N=7$ and $N=11$.  Section \ref{sec:code} at the end of
this paper gives computations of the ramification module for higher
$N$, and the GAP programs we used.

\subsection{The case $N=7$}
\label{example:X(7)a}

The equivalence classes of irreducible representations of $PSL(2,7)$
are $G^*=\{\pi_1,\pi_2,...,\pi_6\}$, where

{\small{
\[
{\rm dim}\, (\pi_1)=1,\ \ \ {\rm dim}\, (\pi_2)={\rm dim}\,
(\pi_3)=3,\ \ \ {\rm dim}\, (\pi_4)=6,\ \ \ {\rm dim}\, (\pi_5)=7,\
\ \ {\rm dim}\, (\pi_6)=8.
\]
}} We can identify
\[
\pi_1=\mathbf{1}_G,\ \ \ \pi_2=X',\ \ \  \pi_3=X'',\ \ \
\pi_4=X_\beta,\ \ \ \pi_5=V,\ \ \ \pi_5=W_\alpha.
\]

Let $\zeta=e^{\frac{2 \pi i}{7}}$ and let $\qqq(q)$ denote the
(quadratic) extension of $\qqq$ by $q=\zeta + \zeta^2 + \zeta^4$.
Let $\mathcal{G}$ denote the Galois group of $\qqq(q)/\qqq$.  Then
$\mathcal{G}$ acts on the irreducible representations $G^*$ by
swapping the two 3-dimensional representations (i.e., $\pi_2=X'$ and
$\pi_3=X''$) and fixing the others.

There are $4$ conjugacy classes of non-trivial cyclic subgroups of
$G$, whose representatives are denoted by $H_1$ (order $2$), $H_2$
(order $3$), $H_3$ (order $7$), $H_4$ (order $4$).  Thus $H_4$ is
spurious. We use GAP to compute the induced characters:

\begin{itemize}
\item
If $\theta_1\in H_1^*$ then
$\pi_{\theta_1}=Ind_{H_1}^G\, \theta_1$ is $84$-dimensional.
Moreover,

\[
\pi_{\theta_1}\cong
\left\{
\begin{array}{ll}
2\pi_2\oplus 2\pi_3\oplus 2\pi_4\oplus 4\pi_5\oplus 4\pi_6, & \theta_1\not= 1,\\
\pi_1\oplus \pi_2\oplus \pi_3\oplus 4\pi_4\oplus 3\pi_5\oplus 4\pi_6, & \theta_1=1.
\end{array}
\right.
\]

\item
If $\theta_2\in H_2^*$ then
$\pi_{\theta_2}=Ind_{H_2}^G\, \theta_2$ is $56$-dimensional.
Moreover,

\[
\pi_{\theta_2}\cong
\left\{
\begin{array}{ll}
\pi_2\oplus \pi_3\oplus 2\pi_4\oplus 2\pi_5\oplus 3\pi_6, & \theta_2\not= 1,\\
\pi_1\oplus \pi_2\oplus \pi_3\oplus 2\pi_4\oplus 3\pi_5\oplus 2\pi_6, & \theta_2=1.
\end{array}
\right.
\]

\item
If $\theta_3\in H_3^*$ is a fixed non-trivial character then
$\pi_{\theta_3}=Ind_{H_3}^G\, \theta_3$ is $24$-dimensional.
Moreover,

\[
\pi_{\theta_3^k}\cong
\left\{
\begin{array}{ll}
\pi_3\oplus \pi_4\oplus \pi_5\oplus \pi_6 & k \text{ quad. non-res. }\pmod 7,\\
\pi_2\oplus \pi_4\oplus \pi_5\oplus \pi_6, & k \text{ quad. res. }\pmod 7,\\
\pi_1 \oplus \pi_5\oplus 2\pi_6, & k\equiv 0\pmod 7.
\end{array}
\right.
\]

\end{itemize}
\noindent These computations agree with the computations in section
\ref{sec:inducedchars}. This data allows us to easily compute the
ramification module using equation (\ref{eqn:ram_mod}):

\begin{equation}
\label{eqn:rammod7}
 [\tilde{\Gamma}_G] =[3\pi_2\oplus 4\pi_3\oplus
6\pi_4\oplus 7\pi_5\oplus 8\pi_6].
\end{equation}

Note that this is not Galois-invariant, because $\pi_2$ and $\pi_3$
have different multiplicities.  A na\"ive computation of the
ramification module, using (\ref{eqn:JKrammod}), yields the following.
For brevity, we represent $m_1[\pi_1]+...+m_6[\pi_6]$ as
$(m_1,...,m_6)$. We compute, using GAP, the quantities

\[
({\rm dim}\, \pi -{\rm dim}\, (\pi^{H_1}))_{i=1..6}
=(1,3,3,6,7,8)-(1,1,1,4,3,4) =(0,2,2,2,4,4),
\]

\[
({\rm dim}\, \pi -{\rm dim}\, (\pi^{H_2}))_{i=1..6}
=(1,3,3,6,7,8)-(1,1,1,2,3,2) =(0,2,2,4,4,6),
\]

\[
({\rm dim}\, \pi -{\rm dim}\, (\pi^{H_3}))_{i=1..6}
=(1,3,3,6,7,8)-(1,0,0,0,1,2) =(0,3,3,6,6,6).
\]

\[
({\rm dim}\, \pi -{\rm dim}\, (\pi^{H_4}))_{i=1..6}
=(1,3,3,6,7,8)-(1,1,1,2,1,2) =(0,2,2,4,6,6).
\]
Combining this with $R_1=R_2=R_3=1$ and $R_4=0$ in
(\ref{eqn:ram_mod}) gives

\[
\begin{array}{c}
[\tilde{\Gamma}_G]= [\bigoplus_{i=1}^6 \left[ \sum_{\ell=1}^4 ({\rm
dim}\, \pi_i -{\rm dim}\, (\pi_i^{H_\ell}))
\frac{R_\ell}{2}\right]\pi_i]\\
=(0,2,2,2,4,4)\frac{1}{2}+(0,2,2,4,4,6)\frac{1}{2}+
(0,3,3,6,6,6)\frac{1}{2}+(0,2,2,4,6,6)\frac{0}{2}\\
=(0,7/2,7/2,6,7,8)\\
=\frac{7}{2}[\pi_2]+ \frac{7}{2}[\pi_3]+ 6[\pi_4]+ 7[\pi_5]+
8[\pi_6].
\end{array}
\]
This is impossible, and therefore we see that $\tilde{\Gamma}_G$
does not have a $\qqq[G]$-module structure in this case.  However,
the result does agree with (\ref{eqn:JKGalclosure}): the
multiplicity given for $\pi_2$ and $\pi_3$ is the average of the two
actual multiplicities.

Now we will use GAP to compute the equivariant degree and
Riemann-Roch module for some example divisors.  As noted in section
\ref{sec:equivdeg}, any effective $G$-invariant divisor on $X(7)$
will be non-special. Since  $X(1)$ is genus zero, for $N=7$,
Borne's formula (\ref{eqn:Borne}) becomes

\[
\begin{array}{ll}
 [L(D)]
&=[\pi_1\oplus 3\pi_2\oplus 3\pi_3\oplus 6\pi_4\oplus 7\pi_5\oplus
  8\pi_6]
+[deg_{eq}(D)]-[\tilde{\Gamma}_G]\\
&=[\pi_1\oplus 3\pi_2\oplus 3\pi_3\oplus 6\pi_4\oplus 7\pi_5\oplus
  8\pi_6] +[deg_{eq}(D)] \\
&\ \ \ \ - 3[\pi_2]- 4[\pi_3]- 6[\pi_4]- 7[\pi_5]- 8[\pi_6] \\
&=[\pi_1] -[\pi_3]+[deg_{eq}(D)].
\end{array}
\]
If $D_1$ is the reduced orbit of a point with stabilizer $H_1$, then
\[
[deg_{eq}(D_1)]=[\pi_{\theta_1}] =[2\pi_2\oplus 2\pi_3\oplus
2\pi_4\oplus 4\pi_5\oplus 4\pi_6]
\]
and
\[
[L(D_1)]=[\pi_1 \oplus 2\pi_2\oplus \pi_3\oplus 2\pi_4\oplus
4\pi_5\oplus 4\pi_6].
\]
If $D_2$ is the reduced orbit of a point with stabilizer $H_2$, then
\[
[deg_{eq}(D_2)]=[\pi_{\theta_2}]= [\pi_2\oplus \pi_3\oplus
2\pi_4\oplus 2\pi_5\oplus 3\pi_6], \]\[ [deg_{eq}(2
D_2)]=[2\pi_{\theta_2}]= [2\pi_2\oplus 2\pi_3\oplus 4\pi_4\oplus
4\pi_5\oplus 6\pi_6],
\]
and
\[
[L(D_2)]=[\pi_1 \oplus \pi_2\oplus 2\pi_4\oplus 2\pi_5\oplus
3\pi_6],
\]\[
[L(2 D_2)]=[\pi_1 \oplus 2\pi_2 \oplus \pi_3 \oplus 4\pi_4 \oplus
4\pi_5\oplus 6\pi_6].
\]
If $D_3$ is the reduced orbit of a point with stabilizer $H_3$, then
\[
[deg_{eq}(D_3)]=[\pi_{\theta_3^{N-1}}]= [\pi_3\oplus \pi_4\oplus
\pi_5\oplus \pi_6], \]\[ [deg_{eq}(2D_3)]=2[\pi_3\oplus \pi_4\oplus
\pi_5\oplus \pi_6],\]\[ [deg_{eq}(3D_3)]=[\pi_2 + 2 \pi_3] +
3[\pi_4\oplus \pi_5\oplus \pi_6],\]\[ [deg_{eq}(4D_3)]=[\pi_2 + 3
\pi_3] + 4[\pi_4\oplus \pi_5\oplus \pi_6], \]\[
[deg_{eq}(5D_3)]=[2\pi_2 + 3 \pi_3] + 5[\pi_4\oplus \pi_5\oplus
\pi_6], \]\[ [deg_{eq}(6D_3)]=[3\pi_2 + 3 \pi_3] + 6[\pi_4\oplus
\pi_5\oplus \pi_6].
\]
It follows that
\[
L(D_3)=[\pi_1]-[\pi_3]+[\pi_3\oplus \pi_4\oplus \pi_5\oplus \pi_6]
=[\pi_1\oplus \pi_4\oplus \pi_5\oplus \pi_6],
\]
which is dimension $22$, and

\[
L(2D_3)=[\pi_1]-[\pi_3]+2[\pi_3\oplus \pi_4\oplus \pi_5\oplus \pi_6]
=[\pi_1\oplus \pi_3]+2[\pi_4\oplus \pi_5\oplus \pi_6],
\]
which is dimension $46$.

\subsection{The case $N=11$}

\label{example:X(11)} We next consider the case of $N=11$. Up to
equivalence, the irreducible representations of $G=PSL(2,11)$ are
$G^*=\{\pi_1,\pi_2,...,\pi_8\}$, where

{\small{
\[
\begin{array}{c}
{\rm dim}\, (\pi_1)=1,\ \ \ {\rm dim}\, (\pi_2)={\rm dim}\,
(\pi_3)=5,\ \ \
{\rm dim}\, (\pi_4)= {\rm dim}\, (\pi_5)=10,\ \ \ \\
{\rm dim}\, (\pi_6)=11,\ \ \ {\rm dim}\, (\pi_7)={\rm dim}\,
(\pi_8)=12.
\end{array}
\]
}} and we can identify
\[ \pi_1=\mathbf{1}_G,\ \ \ \pi_2=X',\ \ \
\pi_3=X'',\ \ \ \pi_4=X_{\beta_1},\ \ \ \ \pi_5=X_{\beta_2},\ \ \ \
\]
\[
\pi_6=V,\ \ \ \pi_7=W_{\alpha_1},\ \ \ \pi_8=W_{\alpha_2}.
\]
We can use GAP to compute the character table of $G$ and from the
character values we can deduce some things about the characters
$\alpha_1, \alpha_2: {\fff}^{\times} \to \ccc^{\times}$ and
$\beta_1, \beta_2: T \to \ccc^{\times}$.  The generator
$\varepsilon$ of ${\fff}^{\times}$ has order 10; it is sent by
$\alpha_1$ to $e^{4 \pi i /5}$ and by $\alpha_2$ to $e^{2 \pi i
/5}$.  The generator $\tau$ of $T$ has order 12; it is sent by
$\beta_1$ to $e^{2 \pi i / 3}$ and by $\beta_2$ to $e^{\pi i/3}$.
Because $N \equiv -1 \pmod {12}$, the numbers $``i"=\tau^{3}$ and
$``\omega"=\tau^2$ are both in $T$; we see from this that
$\beta_1(i)=1$, $\beta_2(i)=-1$, $\beta_1(\omega) \neq 1$, and
$\beta_2(\omega) \neq 1$.

Thanks to the above discussion,
it should be clear that if $K=\qqq(G)$ denotes the
abelian extension of $\qqq$ generated by the character values of $G$
then ${\cal G}=Gal(K/\qqq)$ acts by swapping the two irreducible
$5$-dimensionals $X'$ and $X''$, and the two irreducible
$12$-dimensionals $W_{\alpha_1}$ and $W_{\alpha_2}$, but not the two
irreducible $10$-dimensionals $X_{\beta_1}$ and $X_{\beta_2}$.

There are $5$ conjugacy classes of cyclic subgroups of $G$, whose
representatives are denoted by $H_1$, ..., $H_5$. They satisfy
$|H_1|=2$, $|H_2|=3$, $|H_3|=11$, $|H_4|=5$, $|H_5|=6$.  We are
interested in $H_1$, $H_2$, and $H_3$; $H_4$ and $H_5$ are spurious.
We use GAP to compute the induced characters:

\begin{itemize}
\item
If $\theta_1\in H_1^*$ then $\pi_{\theta_1}=Ind_{H_1}^G\, \theta_1$
is $330$-dimensional. Moreover,

\[
\pi_{\theta_1}\cong \left\{
\begin{array}{ll}
2\pi_2\oplus 2\pi_3\oplus 6\pi_4\oplus 4\pi_5
\oplus 6\pi_6 \oplus 6\pi_7 \oplus 6\pi_8,  & \theta_1\not= 1,\\
\pi_1\oplus 3\pi_2\oplus 3\pi_3\oplus 4\pi_4\oplus 6\pi_5\oplus
5\pi_6 \oplus 6\pi_7\oplus 6\pi_8, & \theta_1=1.
\end{array}
\right.
\]

\item
If $\theta_2\in H_2^*$ then $\pi_{\theta_2} = Ind_{H_2}^G \theta_2$
is $220$-dimensional. Moreover,

\[
\pi_{\theta_2} = \left\{
\begin{array}{ll}
2\pi_2 \oplus 2\pi_3 \oplus 3\pi_4 \oplus 3\pi_5 \oplus
4\pi_6 \oplus 4\pi_7 \oplus 4\pi_8,& \theta_2\not= 1,\\
\pi_1\oplus \pi_2 \oplus \pi_3 \oplus 4\pi_4 \oplus 4\pi_5 \oplus
3\pi_6 \oplus 4\pi_7 \oplus 4\pi_8,& \theta_2= 1.
\end{array}
\right.
\]

\item
If $\theta_3\in H_3^*$ is a fixed non-trivial character
then $\pi_{\theta_3} = Ind_{H_3}^G \theta_3$
is $66$-dimensional. Note that this is not ${\cal G}$-invariant.
Moreover,

\[
\pi_{\theta_3^k}
=
\left\{
\begin{array}{ll}
\pi_2\oplus \pi_4 \oplus \pi_5 \oplus \pi_6 \oplus \pi_7 \oplus \pi_8,
& k \text{ quad. non-res. }\pmod{11},\\
\pi_3\oplus \pi_4 \oplus \pi_5 \oplus \pi_6 \oplus \pi_7 \oplus \pi_8,
& k \text{ quad. res. }\pmod{11}, \\
\pi_1\oplus \pi_6 \oplus 2\pi_7 \oplus 2\pi_8,& k \equiv 0 \pmod{11}.
\end{array}
\right.
\]

\end{itemize}

\noindent From this data, we can easily compute the ramification
module, using equation (\ref{eqn:ram_mod}).  The equivalence class
$[\Gamma_G]$ equals

\[
\begin{array}{c}
\overline{R}_2 [\pi_{\theta_1}]
+\overline{R}_3 ([\pi_{\theta_2}]+2[\pi_{\theta_2^2}])\\
+\overline{R}_N ([\pi_{\theta_3}] +...+ (N-1)[\pi_{\theta_3^{N-1}}])
\end{array}
\]
where $N=11$. As was mentioned above, using \cite{S}, we have
$\overline{R}_2=660/2$, $\overline{R}_3=660/3$,
$\overline{R}_{11}=660/11$. Therefore,

\[
\begin{array}{ll}
[\Gamma_G]&=660\cdot [\frac{1}{2}\cdot (0,2,2,6,4,6,6,6)+
(0,2,2,3,3,4,4,4)\\
&\ \ \ +\cdot (0,3,2,5,5,5,5,5)]\\
&=660\cdot (0,6,5,11,10,12,12,12),
\end{array}
\]
by the formulas for the induced characters above. Therefore,
\newline
$[\tilde{\Gamma}_G]=(0,6,5,11,10,12,12,12)$.

As in the case of $X(7)$, this is not Galois invariant, because
$\pi_2$ and $\pi_3$ have different multiplicities.  Again we compute
the Galois closure using equation (\ref{eqn:JKGalclosure}):

\[
\begin{array}{ll}
({\rm dim}\, \pi -{\rm dim}\, (\pi^{H_1}))_{i=1..8}
&=(1,5,5,10,10,11,12,12)-(1,3,3,4,6,5,6,6)\\
&=(0,2,2,6,4,6,6,6),
\end{array}
\]

\[
\begin{array}{ll}
({\rm dim}\, \pi -{\rm dim}\, (\pi^{H_2}))_{i=1..8}
&=(1,5,5,10,10,11,12,12)-(1,1,1,4,4,3,4,4)\\
&=(0,4,4,6,6,8,8,8),
\end{array}
\]

\[
\begin{array}{ll}
({\rm dim}\, \pi -{\rm dim}\, (\pi^{H_3}))_{i=1..8}
&=(1,5,5,10,10,11,12,12)-(1,0,0,0,0,1,2,2)\\
&=(0,5,5,10,10,10,10,10),
\end{array}
\]

These calculations, combined with $R_1=R_2=R_3=1$ and $R_4=R_5=0$,
give

\[
\begin{array}{c}
[\tilde{\Gamma}_G]= [\bigoplus_{i=1}^8 \left[ \sum_{\ell=1}^5 ({\rm
dim}\, \pi_i -{\rm dim}\, (\pi_i^{H_\ell}))
\frac{R_\ell}{2}\right]\pi_i]\\
=(0,2,2,6,4,6,6,6)\frac{1}{2}+
(0,2,2,3,3,4,4,4)\frac{1}{2}+(0,5,5,10,10,10,10,10)\frac{1}{2}\\
=(0,11/2,11/2,11,10,12,12,12)\\
=\frac{11}{2}[\pi_2]+ \frac{11}{2}[\pi_3]+ 11[\pi_4]+ 10[\pi_5]+ 12[\pi_6]+ 12[\pi_7]+ 12[\pi_8].
\end{array}
\]
We can see directly from this calculation $\tilde{\Gamma}_G$ does
not have a $\qqq[G]$-module structure in this case. However, the
result does agree with (\ref{eqn:JKGalclosure}):  we have computed
the average of the multiplicities of the representations under the
Galois action.

For every non-special $G$-equivariant divisor $D$ on $X(11)$, the
formula (\ref{eqn:Borne}) says

\[
\begin{array}{ll}
 [L(D)]
&=[\pi_1\oplus 5\pi_2\oplus 5\pi_3\oplus 10\pi_4\oplus 10\pi_5
\oplus   11\pi_6\oplus   12\pi_7\oplus   12\pi_8]\\
&\ \ \ \ \ \ \
+[deg_{eq}(D)]-[\tilde{\Gamma}_G]\\
&= [\pi_1]- [\pi_2]-[\pi_4]- [\pi_6]+[deg_{eq}(D)].
\end{array}
\]

If $D_1$ is the reduced orbit of a point with stabilizer $H_1$, then
$D_1$ is non-special and

\[
[deg_{eq}(D_1)]=[\pi_{\theta_1}] =[2\pi_2\oplus 2\pi_3\oplus
6\pi_4\oplus 4\pi_5 \oplus 6\pi_6\oplus 6\pi_7\oplus 6\pi_8],
\]
and
\[
[L(D_1)]=[\pi_1 \oplus \pi_2 \oplus 2\pi_3 \oplus 5\pi_4\oplus
4\pi_5\oplus 5\pi_6\oplus 6\pi_7\oplus 6\pi_8].
\]
If $D_2$ is the reduced orbit of a point with stabilizer $H_2$, then
$D_2$ is non-special and

\[
[deg_{eq}(D_2)]=[\pi_{\theta_2^2}]= [2\pi_2\oplus 2\pi_3\oplus
3\pi_4\oplus 3\pi_5\oplus 4\pi_6 \oplus 4\pi_7\oplus 4\pi_8],
\]
and
\[
[L(D_2)]=[\pi_1 \oplus \pi_2\oplus 2\pi_3\oplus 2\pi_4\oplus
3\pi_5\oplus 3\pi_6\oplus 4\pi_7\oplus 4\pi_8].
\]
If $D_3$ is the reduced orbit of a point with stabilizer $H_3$, then
$D_3$ is non-special and

\[
[deg_{eq}(D_3)]=[\pi_{\theta_3^{10}}]= [\pi_2\oplus \pi_4\oplus
\pi_5\oplus \pi_6\oplus \pi_7\oplus \pi_8],\]\[
[deg_{eq}(2D_3)]=[\pi_2\oplus \pi_3\oplus 2\pi_4\oplus 2\pi_5\oplus
2\pi_6\oplus 2\pi_7\oplus 2\pi_8],\]
and
\[
[L(D_3)]=[\pi_1 \oplus \pi_5\oplus \pi_7 \oplus \pi_8],\]\[
[L(2D_3)]=[\pi_1 \oplus \pi_3\oplus \pi_4\oplus 2\pi_5\oplus
\pi_6\oplus 2\pi_7\oplus 2\pi_8].
\]

\section{Application to codes}

In this section we consider applications of our previous results to
we discuss connections with the theory of error-correcting codes.

Assume that $\ell$ is a good prime. Also, assume that $k$ contains
all the character values of $G$ and that $k$ is finite, where $k$
denotes the field of definition of the reduction of $X$ mod $\ell$.
(The point is that we want to be able to work over a separable
algebraic closure $\overline{k}$ of $k$ but then be able to take
$Gal(\overline{k}/k)$-fixed points to obtain our results.) We recall
some background on AG codes following \cite{JT}.

Let $P_1,...,P_n\in X(k)$ be distinct points
and $E=P_1+...+P_n\in Div(X)$ be stabilized by $G$.
This implies that $G$ acts on the set
$supp(E)$ by permutation.
Assume $D$ is a $G$-equivariant divisor of $X(k)$, so
$G$ acts on the Riemann-Roch space $L(D)$.
Assume these divisors have disjoint support,
${\rm supp}(D)\cap {\rm supp}(E)=\emptyset$.
Let $C=C(D,E)$ denote the AG code

\begin{equation}
\label{eqn:AGcode}
C=\{(f(P_1),...,f(P_n))\ |\ f\in L(D)\}.
\end{equation}
This is the image of $L(D)$ under the evaluation map

\begin{equation}
\label{eqn:eval}
\begin{array}{c}
eval_E:L(D)\rightarrow k^n,\\
f \longmapsto (f(P_1),...,f(P_n)).
\end{array}
\end{equation}
The group $G$ acts on $C$ by $g\in G$ sending
$c=(f(P_1),...,f(P_n))\in C$ to $c'=(f(g^{-1}(P_1)),...,f(g^{-1}(P_n)))$,
where $f\in L(D)$.
First, we observe that this map sending
$c\longmapsto c'$, denoted $\phi(g)$, is well-defined.
In other words, if $eval_E$ is not injective and
$c$ is also represented by $f'\in L(D)$, so
$c=(f'(P_1),...,f'(P_n))\in C$, then we can easily verify
$(f(g^{-1}(P_1)),...,f(g^{-1}(P_n)))=
(f'(g^{-1}(P_1)),...,f'(g^{-1}(P_n)))$.
(Indeed, $G$ acts on the set
$supp(E)$ by permutation.)
This map $\phi(g)$ induces a homomorphism of $G$ into the permutation
automorphism group of the code ${\rm Aut}(C)$, denoted

\begin{equation}
\label{eqn:phi}
\phi:G\rightarrow {\rm Aut}(C)
\end{equation}
For properties of this map, see \cite{JT}.
In particular, the following is known.

\begin{lemma}
If $D$ and $E$ satisfy $\deg(D)>2g$
and $\deg(E)>2g+2$ then $\phi$ and $eval_E$
are injective.
\end{lemma}

\pf
We say that the space $L(D)$
{\bf separates points} if for
all points $P,Q\in X$,
$f(P)=f(Q)$ (for all $f\in L(D)$)
implies $P=Q$ (see \cite{H}, chapter II, \S 7).
By Proposition IV.3.1 in Hartshorne \cite{H},
$D$ very ample implies $L(D)$ separates points.
In general, if $L(D)$ separates points then
\[
{\rm Ker}(\phi)=\{g\in G\ |\ g(P_i)=P_i,\ 1\leq i\leq n\}.
\]
It is known (proof of Prop. VII3.3, \cite{Sti}) that
if $n=\deg(E)>2g+2$ then $\{g\in G\ |\ g(P_i)=P_i,\ 1\leq i\leq n\}$
is trivial. Therefore, if $n>2g+2$ and
$L(D)$ separates points then $\phi$ is injective.
Since (see Corollary IV.3.2 in Hartshorne \cite{H})
$\deg(D)>2g$ implies $D$ is very ample, the lemma follows.
\qed

As an amusing application of our theory, we show how to easily recover some
results of Tsafsman and Vladut on AG codes associated to modular curves.

First, we recall some notation and results from \cite{TV}.
Let $A_N=\zzz [\zeta_N,1/N]$, where $\zeta_N=e^{2\pi i/N}$,
let $K_N$ denote the quadratic subfield of $\qqq (\zeta_N)$,
and let $B_N= A_N\cap K_N$.
There is a scheme $X(N)/\zzz [1/N]$ which represents
a moduli functor ``parameterizing'' elliptic
curves $E$ with a level $N$ structure $\alpha_N$.
There is a scheme $X_P(N)/\zzz [1/N]$ which represents
a moduli functor ``parameterizing'' elliptic
curves $E$ with a ``projective'' level $N$ structure $\beta_N$.
If $P$ is a prime ideal in the ring of integers
${\cal O}_{K_N}$ dividing $\ell$ then the reduction
of the form of $X(N)$ defined over $K_N$, denoted
$X(N)/P$, is a smooth projective absolutely irreducible
curve over the residue field $k(P)$, with a
$PSL_2(\zzz/N\zzz)$-action commuting with the reduction.
Similarly, with $X(N)$ replaced by $X_P(N)$.
Recall from \S 4.1.3 of \cite{TV} that

\[
k(P)=
\left\{
\begin{array}{ll}
GF(\ell^2), & (\ell)=P,\\
GF(\ell), & (\ell)=PP',
\end{array}
\right.
\]
Let $X'_N=X_P(N)/P$ and let

\[
\psi'_N:X'_N\rightarrow X'_N/PSL_2(\zzz/N\zzz)\cong \ppp^1
\]
denote the quotient map.
Let $D_\infty$ denote the reduced orbit (in the sense of Borne)
of the point $\infty$, so $deg(D_\infty)=|G|/N$.
Let $D=rD_\infty$, for $r\geq 1$. According to
\cite{TV}, in general, this divisor is actually defined over
$GF(\ell)$, not just $k(P)$. Moreover,
$\deg(D)=r\cdot (N^2-1)/2$. Let
$E=P_1+...+P_n$ be the sum of all the
supersingular points of $X_N'$ and let

\[
C=C(X'_N,E,D)=\{(f(P_1),...,f(P_n))\ |\
f\in L(D)\}
\]
denote the AG code associated to $X'_N,D,E$.
This is a $G$-module, via (\ref{eqn:phi}).
Moreover, choosing $r$ suitably yields
a ``good'' family of codes with large automorphism group.

In fact, if $D$ is ``sufficiently large'' (so, $D$ is non-special
and both $\phi$, $eval_E$ are injective) then the Brauer-character
analoguess of formulas in section \ref{sec:equivdeg} give not only
the $G$-module structure of each $L(rD_\infty)$, but that of $C$ as
well.

See also Remark 4.1.66 in \cite{TV}.

\section{AG codes associated to $X(7)$}
\label{sec:X7}

We focus on the Klein quartic example started in \cite{JK}.
We also use Elkies \cite{E} as a general reference.

Let $k=GF(43)$. This field contains $7^{th}$ roots of unity
($\zeta_7=41$), cube roots of unity ($\zeta_3=36$), and
the square root of $-7$ (take $\sqrt{-7}=6$).
Consider

\[
\rho_1=\left(
\begin {array}{ccc}
 {\zeta_7}^{4}&0&0\\
 0&{\zeta_7}^{2}&0\\
 0&0&\zeta_7
\end {array}
\right)
=\left(
\begin {array}{ccc}
 16&0&0\\
 0&4&0\\
 0&0&41
\end {array}
\right)
,
\]

\[
\rho_2=\left(
\begin {array}{ccc}
 0&1&0\\
 0&0&1\\
 1&0&0
\end {array}
\right),
\]
and

\[
\begin{array}{ll}
\rho_3&= \left(
\begin {array}{ccc}
 \left( \zeta_7-{\zeta_7}^{6} \right)/\sqrt {-7}&
 \left( {\zeta_7}^{2}-{\zeta_7}^{5} \right)/ \sqrt {-7}&
 \left( {\zeta_7}^{4}-{\zeta_7}^{3} \right)/ \sqrt {-7}\\
 \left( {\zeta_7}^{2}-{\zeta_7}^{5} \right)/ \sqrt {-7}&
 \left( {\zeta_7}^{4}-{\zeta_7}^{3} \right)/ \sqrt {-7}&
 \left( \zeta_7-{\zeta_7}^{6} \right)/ \sqrt {-7}\\
 \left( {\zeta_7}^{4}-{\zeta_7}^{3} \right)/ \sqrt {-7}&
 \left( \zeta_7-{\zeta_7}^{6} \right)/ \sqrt {-7}&
 \left( {\zeta_7}^{2}-{\zeta_7}^{5} \right)/ \sqrt {-7}
\end {array}
\right)\\
&=\left(
\begin {array}{ccc}
11&37&39\\
37&39&11\\
39&11&37
\end {array}
\right).
\end{array}
\]
It may be checked that these matrices preserve the form

\[
\phi(x,y,z)=x^3y+y^3z+z^3x,
\]
over $k$. They generate the subgroup $G\cong PSL_2(7)$
of order $168$ in $PGL(3,k)$. The Klein curve $x^3y+y^3z+z^3x=0$,
denoted here by $X$, has no other automorphisms in
characteristic $43$, so $G=Aut_k(X)$.

Let $D_\infty$ denote the reduced orbit of the point
$\infty$, so $deg(D_\infty)=|G|/N=24$, and let $D=rD_\infty$.
Let $E=P_1+...+P_n$ denote the sum of the
remaining $k(P)$-rational points of $X$,
so $D$ and $E$ have disjoint support.

If $C$ is as in (\ref{eqn:eval}) then
then map $\phi$ in (\ref{eqn:phi}) is injective.
Since $eval_E$ is injective as well, the
$G$-module structure of $C$ is the same as that of
$L(D)$, which is known thanks to the Brauer-character analog
of the formula (\ref{eqn:JKrammod}).
See \S \ref{example:X(7)a} above and also Example 3, \cite{JK}.

Magma tells us that the points of $X(k)$ are

{\footnotesize{
\[
\begin{array}{c}
\{
 (0 : 1 : 0),
 (0 : 0 : 1),
 (1 : 0 : 0),
 (19 : 9 : 1),
 (36 : 9 : 1),
 (31 : 9 : 1),\\
 (19 : 27 : 1),
 (1 : 38 : 1),
 (27 : 38 : 1),
 (15 : 38 : 1),
 (12 : 28 : 1),
 (38 : 28 : 1),\\
 (36 : 28 : 1),
 (40 : 41 : 1),
 (10 : 25 : 1),
 (20 : 25 : 1),
 (13 : 25 : 1),\\
 (20 : 32 : 1),
 (42 : 10 : 1),
 (35 : 10 : 1),
 (9 : 10 : 1),
 (40 : 30 : 1),\\
 (13 : 30 : 1),
 (33 : 30 : 1),
 (24 : 4 : 1),
 (25 : 36 : 1),
 (12 : 36 : 1),
 (6 : 36 : 1),\\
 (12 : 22 : 1),
 (14 : 23 : 1),
 (8 : 23 : 1),
 (21 : 23 : 1),
 (24 : 26 : 1),
 (37 : 26 : 1),\\
 (25 : 26 : 1),
 (23 : 35 : 1),
 (15 : 14 : 1),
 (33 : 14 : 1),
 (38 : 14 : 1),
 (33 : 42 : 1),\\
 (17 : 40 : 1),
 (4 : 40 : 1),
 (22 : 40 : 1),
 (5 : 34 : 1),
 (15 : 34 : 1),
 (23 : 34 : 1),\\
 (31 : 16 : 1),
 (40 : 15 : 1),
 (37 : 15 : 1),
 (9 : 15 : 1),
 (37 : 2 : 1),
 (11 : 6 : 1),\\
 (39 : 6 : 1),
 (36 : 6 : 1),
 (31 : 18 : 1),
 (9 : 18 : 1),
 (3 : 18 : 1),
 (10 : 11 : 1),\\
 (38 : 14 : 1),
 (33 : 42 : 1),\\
 (17 : 40 : 1),
 (4 : 40 : 1),
 (22 : 40 : 1),
 (5 : 34 : 1),
 (15 : 34 : 1),
 (23 : 34 : 1),\\
 (31 : 16 : 1),
 (40 : 15 : 1),
 (37 : 15 : 1),
 (9 : 15 : 1),
 (37 : 2 : 1),
 (11 : 6 : 1),\\
 (39 : 6 : 1),
 (36 : 6 : 1),
 (31 : 18 : 1),
 (9 : 18 : 1),
 (3 : 18 : 1),
 (10 : 11 : 1),\\
 (5 : 13 : 1),
 (24 : 13 : 1),
 (14 : 13 : 1),
 (5 : 39 : 1),
 (41 : 31 : 1),\\
 (13 : 31 : 1),
 (32 : 31 : 1),
 (10 : 7 : 1),
 (14 : 7 : 1),
 (19 : 7 : 1),
 (6 : 21 : 1),\\
 (17 : 17 : 1),
 (3 : 17 : 1),
 (23 : 17 : 1),
 (3 : 8 : 1),
 (25 : 24 : 1),\\
 (16 : 24 : 1),
 (2 : 24 : 1),
 (20 : 29 : 1),
 (17 : 29 : 1),
 (6 : 29 : 1),
 (38 : 1 : 1)
\}
\end{array}
\]
}}
Magma tells us that the orbit of
$(1:0:0)\in X(k)$ under $G$ is

\[
\begin{array}{c}
\{
    ( 1:  0 : 0),
    (36: 28 :1),
    (32: 31 :1),
    (15: 34 :1),
    ( 9: 18 :1),
    (37:  2 :1),\\
    ( 2: 24 :1),
    ( 3:  8 :1),
    (22: 40 :1),
    (19: 27 :1),
    (13: 30 :1),\\
    ( 5: 39 :1),
    ( 8: 23 :1),
    (14:  7 :1),
    (25: 26 :1),
    (33: 42 :1),\\
    (17: 29 :1),
    (42: 10 :1),
    ( 0:  1  0),
    ( 0:  0 :1),\\
    (20: 32 :1),
    (12: 22 :1),
    (27: 38 :1),
    (39:  6 :1)
\}
\end{array}
\]
We evaluate each element of the
Riemann-Roch space $L(D)$ at a point in
the complement of the above-mentioned orbit in
the set of rational points $X(k)$.
For concreteness, a typical basis element of
$L(D)$, when $r=1$, looks like

{\small{
\[
\begin{array}{c}
y^j  x^i  (y + 1)^{-1}  (y + 3)^{-1}  (y + 4)^{-1}  (y + 5)^{-1}  (y + 9)^{-1}
    (y + 11)^{-1}  (y + 12)^{-1}  (y + 13)^{-1}  \times \\
\times    (y + 14)^{-1}  (y + 15)^{-1}  (y + 16)^{-1}  (y + 17)^{-1}  (y +
    19)^{-1}  (y + 20)^{-1}  (y + 21)^{-1} \times \\
\times
(y + 25)^{-1}  (y + 33)^{-1}
(y + 35)^{-1}  (y + 36)^{-1}  (y + 37)^{-1}  (y + 41)^{-1}  \times \\
\times
(x^2y^8 + 34x^2y + 16xy^{11} + 2xy^4 + 17y^{14} + 33y^7 +  17),
\end{array}
\]
}}
for example with $i=2$ and $j=6$. In fact, all such elements
with $i=1,2$ and $0\leq j\leq 6$ are basis elements of $L(D)$.

With $r=1$, Magma says that $C=C(X,D,E)$ is a $[56,22,32]$
code\footnote{In other words, $C$ has length $56$,
dimension $22$ over $k$, and minimum distance $32$.}.
It has a generator matrix of in standard form
$(I\ |\ A)$, where $I$ is the $22\times 22$ identity and
$A$ can be given explicitly as well.
In this case, $eval_E$ is injective.

With $r=2$, Magma says that $C=C(X,D,E)$ is a $[56,46,8]$ code.
It has a generator matrix of in standard form
$(I\ |\ A)$, where $I$ is the $46\times 46$ identity and
$A$ can be given explicitly as well.
In this case, $eval_E$ is injective.

With $r=3$, Magma says that $C=C(X,D,E)$ is a $[56,56,1]$ code.
In this case, $eval_E$ is not injective. Indeed, $\dim L(3D_\infty)=70$.

\begin{remark}
Indeed, it is known more generally, that for an AG code
constructed as above from a curve of genus $g$ that
$n\leq \dim(C)+d(C)+g-1$, where $d(C)$ denotes
the minimum distance (Theorem 3.1.1 in \cite{TV}).
Therefore, as an AG code, the codes constructed above with
$r=1,2$ are as ``long as possible''.
Are they best possible over $GF(43)$ for those
$k,d$?
\end{remark}

\begin{remark}
An attempt to compute the permutation automorphism group of the
above $[56,22,32]$ code using MAGMA failed due to lack
of memory (on an AMD64 linux machine with 1.5 G of RAM).
\end{remark}

\begin{remark}
In general, \S 4.1 of \cite{TV}
shows how to construct a family of `good'' codes from
the curves $X=X_N'$, for $N>5$ is a prime, with
automorphism group $G=PSL(2,p)$.
\end{remark}

\section{Appendix: Tables for Theorem \ref{thrm:main}}
\label{sec:appendix}

The following tables are an intermediate step in the proof of
theorem \ref{thrm:main}. They were obtained by compiling the results
in \S\S \ref{sec:IndH1}-\ref{sec:IndH3}.

\[
\begin{array}[ht]{|c|c|c| c|}\hline
&&& \\
N\equiv 1\pmod{24} & Ind_{H_1}^G\theta_1 &
Ind_{H_2}^G(\theta_2+2\theta_2^2) &
Ind_{H_3}^G(\sum_{\ell =1}^{N-1}\ell \theta_3^\ell)\\
&&& \\ \hline
\mathbf{1}_G & 0 & 0 & 0 \\
W' & \frac{N-1}{4} & \frac{N-1}{2} & S_{\cal Q} \\
W'' & \frac{N-1}{4} & \frac{N-1}{2} & S_{\cal N} \\
\begin{array}{c}
W_{\alpha} {\rm \ such\ that}\\
\alpha(i)=\alpha(\omega)=1
\end{array} & \frac{N-1}{2} & N-1  & \frac{N(N-1)}{2} \\
\begin{array}{c}
W_{\alpha} {\rm \ such\ that}\\
\alpha(i)=1,\alpha(\omega)\not =1
\end{array} & \frac{N-1}{2} & N+2  & \frac{N(N-1)}{2} \\
\begin{array}{c}
W_{\alpha} {\rm \ such\ that}\\
\alpha(i)=-1,\alpha(\omega)=1
\end{array} & \frac{N+3}{2} & N-1  & \frac{N(N-1)}{2} \\
\begin{array}{c}
W_{\alpha} {\rm \ such\ that}\\
\alpha(i)=-1,\alpha(\omega)\not= 1
\end{array} & \frac{N+3}{2} & N+2 & \frac{N(N-1)}{2} \\
X_{\beta}  & \frac{N-1}{2} & N-1  & \frac{N(N-1)}{2} \\
V & \frac{N-1}{2} & N-1 & \frac{N(N-1)}{2} \\ \hline
\end{array}
\]

\[
\begin{array}[ht]{|c|c|c| c|}\hline
&&& \\
N\equiv 5\pmod{24} & Ind_{H_1}^G\theta_1 &
Ind_{H_2}^G(\theta_2+2\theta_2^2) &
Ind_{H_3}^G(\sum_{\ell =1}^{N-1}\ell \theta_3^\ell)\\
&&& \\ \hline
\mathbf{1}_G & 0 & 0 & 0 \\
W' & \frac{N+3}{4} & \frac{N+1}{2} & S_{\cal Q} \\
W'' & \frac{N+3}{4} & \frac{N+1}{2} & S_{\cal N} \\
\begin{array}{c}
W_{\alpha} {\rm \ such\ that}\\
\alpha(i)=1
\end{array} & \frac{N-1}{2} & N+1  & \frac{N(N-1)}{2} \\
\begin{array}{c}
W_{\alpha} {\rm \ such\ that}\\
\alpha(i)=-1
\end{array} & \frac{N+3}{2} & N+1  & \frac{N(N-1)}{2} \\
\begin{array}{c}
X_{\beta} {\rm \ such\ that}\\
\beta(\omega) =1
\end{array} & \frac{N-1}{2} & N+1  & \frac{N(N-1)}{2} \\
\begin{array}{c}
X_{\beta} {\rm \ such\ that}\\
\beta(\omega)\not =1
\end{array} & \frac{N-1}{2} & N-2  & \frac{N(N-1)}{2} \\
V & \frac{N-1}{2} & N+1 & \frac{N(N-1)}{2} \\ \hline
\end{array}
\]

\[
\begin{array}[ht]{|c|c|c| c|}\hline
&&& \\
N\equiv 7\pmod{24} & Ind_{H_1}^G\theta_1 &
Ind_{H_2}^G(\theta_2+2\theta_2^2) &
Ind_{H_3}^G(\sum_{\ell =1}^{N-1}\ell \theta_3^\ell)\\
&&& \\ \hline
\mathbf{1}_G & 0 & 0 & 0 \\
\begin{array}{c}
W_{\alpha} {\rm \ such\ that}\\
\alpha(\omega)=1
\end{array} & \frac{N+1}{2} & N-1  & \frac{N(N-1)}{2} \\
\begin{array}{c}
W_{\alpha} {\rm \ such\ that}\\
\alpha(\omega)\not =1
\end{array} & \frac{N+1}{2} & N+2  & \frac{N(N-1)}{2} \\
X'  & \frac{N+1}{4} & \frac{N-1}{2} & S_{\cal Q} \\
X'' & \frac{N+1}{4} & \frac{N-1}{2} & S_{\cal N} \\
\begin{array}{c}
X_{\beta} {\rm \ such\ that}\\
\beta(i)=1
\end{array} & \frac{N+1}{2} & N-1  & \frac{N(N-1)}{2} \\
\begin{array}{c}
X_{\beta} {\rm \ such\ that}\\
\beta(i)=-1
\end{array} & \frac{N-3}{2} & N-1  & \frac{N(N-1)}{2} \\
V & \frac{N+1}{2} & N-1 & \frac{N(N-1)}{2} \\ \hline
\end{array}
\]

\[
\begin{array}[ht]{|c|c|c| c|}\hline
&&& \\
N\equiv 11\pmod{24} & Ind_{H_1}^G\theta_1 &
Ind_{H_2}^G(\theta_2+2\theta_2^2) &
Ind_{H_3}^G(\sum_{\ell =1}^{N-1}\ell \theta_3^\ell)\\
&&& \\ \hline
\mathbf{1}_G & 0 & 0 & 0 \\
W_{\alpha} & \frac{N+1}{2} & N+1  & \frac{N(N-1)}{2} \\
X'  & \frac{N-3}{4} & \frac{N+1}{2} & S_{\cal Q} \\
X'' & \frac{N-3}{4} & \frac{N+1}{2} & S_{\cal N} \\

\begin{array}{c}
X_{\beta} {\rm \ such\ that}\\
\beta(i)=\beta(\omega) =1
\end{array} & \frac{N+1}{2} & N+1  & \frac{N(N-1)}{2} \\
\begin{array}{c}
X_{\beta} {\rm \ such\ that}\\
\beta(i)=1,\beta(\omega)\not =1
\end{array} & \frac{N+1}{2} & N-2  & \frac{N(N-1)}{2} \\
\begin{array}{c}
X_{\beta} {\rm \ such\ that}\\
\beta(i)=-1,\beta(\omega) =1
\end{array} & \frac{N-3}{2} & N+1  & \frac{N(N-1)}{2} \\
\begin{array}{c}
X_{\beta} {\rm \ such\ that}\\
\beta(i)=-1,\beta(\omega)\not =1
\end{array} & \frac{N-3}{2} & N-2  & \frac{N(N-1)}{2} \\

V & \frac{N+1}{2} & N+1 & \frac{N(N-1)}{2} \\ \hline
\end{array}
\]

\[
\begin{array}[ht]{|c|c|c| c|}\hline
&&& \\
N\equiv 13\pmod{24} & Ind_{H_1}^G\theta_1 &
Ind_{H_2}^G(\theta_2+2\theta_2^2) &
Ind_{H_3}^G(\sum_{\ell =1}^{N-1}\ell \theta_3^\ell)\\
&&& \\ \hline
\mathbf{1}_G & 0 & 0 & 0 \\
W' & \frac{N+3}{4} & \frac{N-1}{2} & S_{\cal Q} \\
W'' & \frac{N+3}{4} & \frac{N-1}{2} & S_{\cal N} \\
\begin{array}{c}
W_{\alpha} {\rm \ such\ that}\\
\alpha(i)=\alpha(\omega)=1
\end{array} & \frac{N-1}{2} & N-1  & \frac{N(N-1)}{2} \\
\begin{array}{c}
W_{\alpha} {\rm \ such\ that}\\
\alpha(i)=1,\alpha(\omega)\not =1
\end{array} & \frac{N-1}{2} & N+2  & \frac{N(N-1)}{2} \\
\begin{array}{c}
W_{\alpha} {\rm \ such\ that}\\
\alpha(i)=-1,\alpha(\omega)=1
\end{array} & \frac{N+3}{2} & N-1  & \frac{N(N-1)}{2} \\
\begin{array}{c}
W_{\alpha} {\rm \ such\ that}\\
\alpha(i)=-1,\alpha(\omega)\not= 1
\end{array} & \frac{N+3}{2} & N+2 & \frac{N(N-1)}{2} \\
X_{\beta}  & \frac{N-1}{2} & N-1  & \frac{N(N-1)}{2} \\
V & \frac{N-1}{2} & N-1 & \frac{N(N-1)}{2} \\ \hline
\end{array}
\]

\[
\begin{array}[ht]{|c|c|c| c|}\hline
&&& \\
N\equiv 17\pmod{24} & Ind_{H_1}^G\theta_1 &
Ind_{H_2}^G(\theta_2+2\theta_2^2) &
Ind_{H_3}^G(\sum_{\ell =1}^{N-1}\ell \theta_3^\ell)\\
&&& \\ \hline
\mathbf{1}_G & 0 & 0 & 0 \\
W'  & \frac{N-1}{4} & \frac{N+1}{2} & S_{\cal Q} \\
W'' & \frac{N-1}{4} & \frac{N+1}{2} & S_{\cal N} \\
\begin{array}{c}
W_{\alpha} {\rm \ such\ that}\\
\alpha(i)=1
\end{array} & \frac{N-1}{2} & N+1  & \frac{N(N-1)}{2} \\
\begin{array}{c}
W_{\alpha} {\rm \ such\ that}\\
\alpha(i)=-1
\end{array} & \frac{N+3}{2} & N+1  & \frac{N(N-1)}{2} \\
\begin{array}{c}
X_{\beta} {\rm \ such\ that}\\
\beta(\omega) =1
\end{array} & \frac{N-1}{2} & N+1  & \frac{N(N-1)}{2} \\
\begin{array}{c}
X_{\beta} {\rm \ such\ that}\\
\beta(\omega)\not =1
\end{array} & \frac{N-1}{2} & N-2  & \frac{N(N-1)}{2} \\
V & \frac{N-1}{2} & N+1 & \frac{N(N-1)}{2} \\ \hline
\end{array}
\]

\[
\begin{array}[ht]{|c|c|c| c|}\hline
&&& \\
N\equiv 19\pmod{24} & Ind_{H_1}^G\theta_1 &
Ind_{H_2}^G(\theta_2+2\theta_2^2) &
Ind_{H_3}^G(\sum_{\ell =1}^{N-1}\ell \theta_3^\ell)\\
&&& \\ \hline
\mathbf{1}_G & 0 & 0 & 0 \\
\begin{array}{c}
W_{\alpha} {\rm \ such\ that}\\
\alpha(\omega)=1
\end{array} & \frac{N+1}{2} & N-1  & \frac{N(N-1)}{2} \\
\begin{array}{c}
W_{\alpha} {\rm \ such\ that}\\
\alpha(\omega)\not =1
\end{array} & \frac{N+1}{2} & N+2  & \frac{N(N-1)}{2} \\
X'  & \frac{N-3}{4} & \frac{N-1}{2} & S_{\cal Q} \\
X'' & \frac{N-3}{4} & \frac{N-1}{2} & S_{\cal N} \\

\begin{array}{c}
X_{\beta} {\rm \ such\ that}\\
\beta(i)=1
\end{array} & \frac{N+1}{2} & N-1  & \frac{N(N-1)}{2} \\
\begin{array}{c}
X_{\beta} {\rm \ such\ that}\\
\beta(i)=-1
\end{array} & \frac{N-3}{2} & N-1  & \frac{N(N-1)}{2} \\
V & \frac{N+1}{2} & N-1 & \frac{N(N-1)}{2} \\ \hline
\end{array}
\]

\[
\begin{array}[ht]{|c|c|c| c|}\hline
&&& \\
N\equiv 23\pmod{24} & Ind_{H_1}^G\theta_1 &
Ind_{H_2}^G(\theta_2+2\theta_2^2) &
Ind_{H_3}^G(\sum_{\ell =1}^{N-1}\ell \theta_3^\ell)\\
&&& \\ \hline
\mathbf{1}_G & 0 & 0 & 0 \\
W_{\alpha}  & \frac{N+1}{2} & N+1  & \frac{N(N-1)}{2} \\
X'  & \frac{N+1}{4} & \frac{N+1}{2} & S_{\cal Q} \\
X'' & \frac{N+1}{4} & \frac{N+1}{2} & S_{\cal N} \\

\begin{array}{c}
X_{\beta} {\rm \ such\ that}\\
\beta(i)=\beta(\omega) =1
\end{array} & \frac{N+1}{2} & N+1  & \frac{N(N-1)}{2} \\
\begin{array}{c}
X_{\beta} {\rm \ such\ that}\\
\beta(i)=1,\beta(\omega)\not =1
\end{array} & \frac{N+1}{2} & N-2  & \frac{N(N-1)}{2} \\
\begin{array}{c}
X_{\beta} {\rm \ such\ that}\\
\beta(i)=-1,\beta(\omega) =1
\end{array} & \frac{N-3}{2} & N+1  & \frac{N(N-1)}{2} \\
\begin{array}{c}
X_{\beta} {\rm \ such\ that}\\
\beta(i)=-1,\beta(\omega)\not =1
\end{array} & \frac{N-3}{2} & N-2  & \frac{N(N-1)}{2} \\
V & \frac{N+1}{2} & N+1 & \frac{N(N-1)}{2} \\ \hline
\end{array}
\]

\section{Appendix: Ramification modules and GAP code}
\label{sec:code}

We used GAP to compute the ramification module for many small values
of $N$.  Here are the first results:

{\footnotesize{
\begin{center}
\begin{tabular}{|c|c|c|}\hline
 & & \\
$N$ & $G$-module structure of $\tilde{\Gamma_G}$ & $G$-module structure of $\tilde{\Gamma_G}$\\
 &   by definition & using JK formula \\ \hline
5 & [ 0, 3, 3, 4, 5 ] & [ 0, 3, 3, 4, 5 ] \\
7 & [ 0, 3, 4, 6, 7, 8 ] & [ 0, 7/2, 7/2, 6, 7, 8 ] \\
11 & [ 0, 5, 6, 11, 10, 12, 12, 12 ] &[ 0, 11/2, 11/2, 11, 10, 12, 12, 12 ]  \\
13 & [ 0, 7, 7, 13, 13, 13, 13, 14, 15 ]
& [ 0, 7, 7, 13, 13, 13, 13, 14, 15 ] \\
17 &[ 0, 9, 9, 18, 17, 17, 17, 18,  ... ]
& [ 0, 9, 9, 18, 17, 17, 17, 18,  ... ]\\
19 &[ 0, 9, 10, 20, 20, 19, 19, 20, ... ]
& [ 0, 19/2, 19/2, 20, 20, 19, 19, 20,  ... ] \\
23 & [ 0, 11, 14, 24, 24, 24, 23, 23, ... ]
& [ 0, 25/2, 25/2, 24, 24, 24, 23, 23, ... ]\\
29 &
[ 0, 16, 16, 30, 31, 31, 30, ... ]
&[ 0, 16, 16, 30, 31, 31, 30,  ... ]\\
\hline
\end{tabular}
\end{center}
}}

GAP was able to go up to and including $N=277$
(in which case $G$ is order about 20 million).

\vskip .3in

GAP code for the tables

{\footnotesize{
\begin{verbatim}


ram_module_X:=function(p)
## p is a prime
## output is (m1,...,mn)
## where n = # conj classes of G=PSL(2,p)
## and mi = mult of pi_i in ram mod of G
##
local G,i,j,n,H,G1,H_chars,CG,G_chars,w,m,theta,pi_theta_to_the;
  G:=PSL(2,p);
  H:=[];
  H_chars:=[];
  CG:=ConjugacyClassesSubgroups(G);
  H[1]:=Representative(CG[2]); # size 2
  H[2]:=Representative(CG[3]); # size 3
  n :=Size(CG);
  for i in [1..n] do
    if Size(Representative(CG[i]))=p then
      H[3]:=Representative(CG[i]);
    fi;
  od;
  for i in [1..Size(H)] do
   H_chars[i]:=Irr(H[i]);
  od;
  G_chars:=Irr(G);
  m:=[];
  m[1]:=[];m[2]:=[];m[3]:=[];
  theta:=List([1..3],i->H_chars[i][2]);
  pi_theta_to_the:=List([1..3],i->Sum([1..(Size(H[i])-1)],
      j->j*InducedClassFunction(theta[i]^(j),G)));
  for i in [1..3] do
    m[i]:=List(G_chars, pi->ScalarProduct(pi_theta_to_the[i],pi))/Size(H[i]);
  od;
Print("\n\n m = ",m,"\n\n");
return Sum(m);
end;

ram_module_X_JK:=function(p)
## p is a prime
## output is (m1,...,mn)
## where n = # conj classes of G=PSL(2,p)
## and mi = "mult of pi_i in ram mod of G" using JK formula
##
local G,i,j,n,H,G1,H_chars,CG,G_chars,A,B,C,D,pi;
  G:=PSL(2,p);
  H:=[];
  H_chars:=[];
  CG:=ConjugacyClassesSubgroups(G);
  H[1]:=Representative(CG[2]); # size 2
  H[2]:=Representative(CG[3]); # size 3
  n:=Size(CG);
  for i in [1..n] do
    if Size(Representative(CG[i]))=p then
      H[3]:=Representative(CG[i]);
    fi;
  od;
  for i in [1..Size(H)] do
   H_chars[i]:=Irr(H[i]);
  od;
  G_chars:=Irr(G);
  n:=Length(G_chars);
  A:=[];
  for i in [1..n] do
    pi:=G_chars[i];
    A[i]:=ScalarProduct(H_chars[1][1],RestrictedClassFunction(pi,H[1]));
  od;
  B:=[];
  for i in [1..n] do
    pi:=G_chars[i];
    B[i]:=ScalarProduct(H_chars[2][1],RestrictedClassFunction(pi,H[2]));
  od;
  C:=[];
  for i in [1..n] do
    pi:=G_chars[i];
    C[i]:=ScalarProduct(H_chars[3][1],RestrictedClassFunction(pi,H[3]));
  od;
  D:=[];
  for i in [1..n] do
    pi:=G_chars[i];
    D[i]:=DegreeOfCharacter(pi);
  od;
 return (1/2)*(3*D-A-B-C);
end;


equiv_deg_module_X:=function(p,ii,r)
## p is a prime
## ii  = 1 for H[1], size 2
## ii  = 2 for H[2], size 3
## ii  = 3 for H[3], size p
## output is (m1,...,mn)
## where n = # conj classes of G=PSL(2,p)
## and mi = mult of pi_i in deg_equiv module of G
##
local G,i,j,n,H,G1,H_chars,CG,G_chars,w,m,theta,pi_theta,pi_theta_to_the;
G:=PSL(2,p);
H:=[]; # 3 decomp gps
H_chars:=[];
CG:=ConjugacyClassesSubgroups(G);
H[1]:=Representative(CG[2]); # size 2
H[2]:=Representative(CG[3]); # size 3
n:=Size(CG);
for i in [1..n] do
  if Size(Representative(CG[i]))=p then
    H[3]:=Representative(CG[i]);
  fi;
od;
for i in [1..3] do
 H_chars[i]:=Irr(H[i]);
od;
G_chars:=Irr(G);
m:=[];
m[1]:=[];m[2]:=[];m[3]:=[];
theta:=List([1..3],i->H_chars[i][2]);
pi_theta_to_the:=List([1..3],i->List([1..r],j->InducedClassFunction(theta[i]^(-j),G)));
for i in [1..3] do
 for j in [1..r] do
  m[i][j]:=List(G_chars, pi->ScalarProduct(pi_theta_to_the[i][j],pi));
 od;
od;
return Sum(m[ii]);
end;

\end{verbatim}
}}

\end{document}